\input ./style/arxiv-general.cfg
\documentclass[MSNbibl,number,citesort,seceqn,dvips]{arxbj}
\makeatletter
   \@ifpackageloaded{graphicx}{}{\usepackage{graphicx}}
\makeatother
\usepackage{upgreek,dcolumn}

\volume{22}
\issue{3}
\pubyear{2016}
\firstpage{1301}
\lastpage{1330}
\doi{10.3150/14-BEJ692}
\docsubty{FLA}

\makeatletter
\newcolumntype{d}[1]{D{.}{.}{#1}}
\def\pi{\uppi}
\newcommand{\rrvert}{\vert}
\newcommand{\rrVert}{\Vert}
\newcommand{\llvert}{\vert}
\newcommand{\llVert}{\Vert}
\newcommand{\eqref}[1]{(\ref{#1})}
\newtheorem{thmm}{Theorem}[section]
\newremark{rem}{Remark}[section]
\newtheorem{lemma}{Lemma}[section]
\newtheorem{prop}{Proposition}[section]
\newtheorem{cor}{Corollary}[section]
\makeatother

\begin{document}
\begin{frontmatter}

\title{Estimation of inverse autocovariance matrices for long memory processes}
\runtitle{Estimation of inverse autocovariance matrices}

\begin{aug}
\author[A]{\inits{C.-K.}\fnms{Ching-Kang} \snm{Ing}\thanksref{A,e1}\ead[label=e1,mark]{cking@stat.sinica.edu.tw}},
\author[B]{\inits{H.-T.}\fnms{Hai-Tang} \snm{Chiou}\thanksref{B,e2}\ead[label=e2,mark]{d002040001@student.nsysu.edu.tw}}
\and
\author[B]{\inits{M.}\fnms{Meihui} \snm{Guo}\thanksref{B,e3}\ead[label=e3,mark]{guomh@math.nsysu.edu.tw}}

\address[A]{Institute of Statistical Science,
Academia Sinica, Taipei, Taiwan.\\
\printead{e1}}

\address[B]{Department of Applied Mathematics,
National Sun Yat-sen University, Kaohsiung, Taiwan.\\
\printead{e2,e3}}


\end{aug}

%
\received{\smonth{12} \syear{2013}}
%
\revised{\smonth{11} \syear{2014}}

%
\begin{abstract}
This work aims at estimating inverse autocovariance matrices of long
memory processes admitting a linear representation. A modified Cholesky
decomposition is used in
conjunction with an increasing order autoregressive model to achieve
this goal. The spectral norm consistency of the proposed estimate is
established. We then extend this result to linear regression models with
long-memory time series errors. In particular, we show that when the
objective is to consistently estimate the inverse autocovariance
matrix of the error process, the same approach still works well if
the estimated (by least squares) errors are used in place of the
unobservable ones.
Applications of this result
to estimating unknown parameters in the aforementioned regression model
are also given.
Finally, a simulation study is performed to
illustrate our theoretical findings.
\end{abstract}

%
\end{frontmatter}

\section{Introduction}\label{section1}
Statistical inference for dependent data
often involves consistent estimates of
the inverse autocovariance matrix
of a stationary time series.
For example,
by making use of
a consistent estimate of
the inverse autocovariance matrix
of a short-memory time series
(in the sense that its autocovariance function is absolutely summable),
Wu and Pourahmadi \cite{r17} constructed estimates
of the finite-past predictor coefficients of the time series and
derived their error bounds.
Moreover, in regression models with short-memory errors,
Cheng, Ing and Yu \cite{r1114}
proposed feasible generalized least squares estimates (FGLSE)
of the regression coefficients
using a consistent estimate of the inverse autocovariance matrix
of the error process.
They then established an asymptotically efficient
model averaging result based on the FGLSEs.

Having observed a realization $u_{1}, \ldots, u_{n}$
of a zero-mean stationary time series $\{u_{t}\}$,
a natural estimate of its autocovariance function $\gamma_{k}=\operatorname
{cov}(u_{0}, u_{k})$
is the sample autocovariance function $\hat{\gamma}_{k}=n^{-1}\sum_{i=1}^{n-|k|}u_{i}u_{i+|k|}$, $k=0, \pm1, \ldots, \pm(n-1)$.
Moreover, it is known that the $k_{n}$-dimensional sample
autocovariance matrix
$\breve{\bolds{\Omega}}_{k_{n}}=(\hat{\gamma}_{i-j})_{1\leq i,
j\leq k_{n}}$
and its inverse $\breve{\bolds{\Omega}}{}^{-1}_{k_{n}}$
are consistent estimates of their population counterparts
$\bolds{\Omega}_{k_{n}}=(\gamma_{i-j})_{1\leq i, j\leq k_{n}}$
and $\bolds{\Omega}^{-1}_{k_{n}}$,
provided $k_{n}\ll n$ and $\sum_{k=1}^{\infty}|\gamma_{k}|< \infty$.
See, for example, Berk \cite{r2}, Shibata \cite{r14}, Ing and Wei \cite
{r9} and Wu and Pourahmadi \cite{r17}.
However, when the objective is to estimate the
$n$-dimensional autocovariance matrix $\bolds{\Omega}_{n}$,
Wu and Pourahmadi \cite{r17} showed that $\breve{\bolds{\Omega}}_{n}$
is no longer consistent in the short-memory case.
In addition, Palma and Pourahmadi \cite{r13} pointed out that this dilemma
carries over to the long-memory case, assuming $\sum_{k=1}^{\infty
}|\gamma_{k}|=\infty$.
To circumvent this
difficulty,
Wu and Pourahmadi \cite{r17} proposed a banded covariance matrix estimate
$\breve{\bolds{\Omega}}_{n, l}=(\hat{\gamma}_{i-j}\mathbf
{1}_{|i-j|\leq l})_{1\leq i, j\leq n}$
of $\bolds{\Omega}_{n}$, where $l\geq0$ is an integer
and called the banding parameter.
When $\{u_{t}\}$ is a short-memory time series satisfying some mild conditions
and $l=l_{n}$ grows to infinity with $n$ at a suitable rate, they
established consistency of
$\breve{\bolds{\Omega}}_{n, l}$ and $\breve{\bolds{\Omega
}}{}^{-1}_{n, l}$
under spectral norm.
The result of Wu and Pourahmadi \cite{r17} was subsequently
improved by Xiao and Wu \cite{r0} to a better convergence rate,
and extended by
McMurry and Politis \cite{r12} to tapered covariance matrix estimates.
Alternatively,
Bickel and Gel~\cite{r31} considered a banded covariance matrix
estimate $\breve{\bolds{\Omega}}_{p_{n}, l}$
of $\bolds{\Omega}_{p_{n}}$, with $p_{n}=\mathrm{o}(n)$.
Assuming that $\{u_{t}\}$ is a stationary short-memory
AR($\infty$) process,
they obtained $\breve{\bolds{\Omega}}_{p_{n}, l}$'s consistency under
the Frobenius norm, provided $l=l_{n}$ tends to infinity sufficiently slowly.

Although the banded and tapered covariance matrix estimates work well
for the short-memory time series,
they are not necessarily suitable for the long-memory case because the
autocovariance function of the latter
is not absolutely summable.
As a result, the banded and tapered matrix estimates
may incur large truncation errors,
which prevent them from achieving consistency.
A~major repercussion of this inconsistency property
is that a consistent estimate of $\bolds{\Omega}^{-1}_{n}$ can no
longer be obtained
by inverting $\breve{\bolds{\Omega}}_{n, l}$ or its tapered version.
On the other hand, since the spectral densities
of most long-memory time series encountered in common practice
are bounded away from zero, it follows from Proposition~4.5.3 of
Brockwell and Davis \cite{r4}
that
%
\begin{equation}
\sup_{k\geq
1}\bigl\|\bolds{\Omega}_{k}^{-1}
\bigr\|_{2}< \infty\label{lemma3.2-3},
\end{equation}
where for a $k$-dimensional matrix $A$,
$\|A\|_{2}=\sup_{\{\mathbf{x}\in R^{k}\dvt \mathbf{x}'\mathbf{x}=1\}
}(\mathbf{x}'\mathbf{A}'\mathbf{Ax})^{1/2}$ denotes its spectral norm.
Motivated by \eqref{lemma3.2-3},
this paper aims to propose
a direct estimate of
$\bolds{\Omega}^{-1}_{n}$ and establish its consistency in the
spectral norm sense,
which is particularly relevant under the long-memory setup.

To fix ideas, assume
%
\begin{equation}
u_{t}=\sum_{j=0}^{\infty}
\psi_{j}w_{t-j}, \label{model}
\end{equation}
where $\psi_{0}=1$ and $\{w_{t}\}$ is a martingale difference sequence
with $E(w_{t})=0$ and $E(w^{2}_{t})=\sigma^{2}$ for all $t$,
and
%
\begin{equation}
\psi_{j}=\mathrm{O}\bigl(j^{-1+d}\bigr), \label{model-1}
\end{equation}
with $d$ satisfying $0<d<1/2$.
We shall also assume that
$\{u_{t}\}$ admits
an AR$(\infty)$ representation,
%
\begin{equation}
u_{t}=\sum_{i=1}^{\infty}a_{i}u_{t-i}+w_{t},\label{model2}
\end{equation}
where
%
\begin{equation}
a_{i}=\mathrm{O}\bigl(i^{-1-d}\bigr).\label{model-2}
\end{equation}
In view of \eqref{model-1},
the autocovariance function of $\{u_{t}\}$ obeys
%
\begin{equation}
\gamma_{k}=\sum_{j=0}^{\infty}
\psi_{j}\psi_{j+|k|}\sigma ^{2}=\mathrm{O}
\bigl(|k|^{-1+2d}\bigr), \label{acf}
\end{equation}
which may not be absolutely summable.
A well-known model satisfying \eqref{model}--\eqref{model-2} is the
FARIMA$(p,d,q)$ processes,
%
\begin{equation}
\phi(B) (1-B)^d u_t=\theta(B) w_t,
\label{farima}
\end{equation}
where $B$ is the backshift operator, $\phi(z)$ and $\theta(z)$ are
polynomials of orders $p$ and $q$, respectively,
$|\phi(z)\theta(z)|\neq0$ for $|z| \leq1$, and $|\phi(z)|$ and
$|\theta(z)|$
have no common zeros.
Note that when \eqref{farima} is assumed, the spectral density of $\{
u_t\}$, $f_{u}(\lambda)$,
satisfies
%
\begin{equation}
\inf_{\lambda\in[-\pi, \pi]}f_{u}(\lambda) > 0, \label{model3}
\end{equation}
from which
\eqref{lemma3.2-3} follows.

Let
%
\begin{equation}
\sigma_{k}^{2}=E(u_{t}-a_{k,1}u_{t-1}-
\cdots-a_{k,k}u_{t-k})^{2}, \label{sigmak}
\end{equation}
where
$k \geq1$ and
%
\begin{equation}
(a_{k,1},\ldots,a_{k,k})=\mathop{\operatorname{arg\,min}}_{(\alpha_{1},\ldots,\alpha
_{k})\in
R^{k}}E(u_{t}-
\alpha_{1}u_{t-1}-\cdots-\alpha_{k}u_{t-k})^{2}
\label{regcoe}.
\end{equation}
To directly estimate $\bolds{\Omega}^{-1}_{n}$, we start by
defining the modified Cholesky decomposition (see, e.g., Berk \cite{r2}
and Wu and Pourahmadi \cite{r16})
of $\bolds{\Omega}_{n}$:
\[
\mathbf{T}_{n}\bolds{\Omega}_{n}\mathbf{T}_{n}'=
\mathbf {D}_{n}\label{autoCov},
\]
where
\[
\mathbf{D}_{n}=\operatorname{diag}\bigl(\gamma_{0},
\sigma_{1}^{2}, \sigma_{2}^{2},
\ldots, \sigma_{n-1}^{2}\bigr),
\]
and
$\mathbf{T}_{n}=(t_{ij})_{1\leq i,j\leq n}$ is a lower triangular matrix
satisfying
\[
t_{ij}= \cases{ 0,&\quad$\mbox{if } i < j$;\vspace*{2pt}
\cr
1,&\quad$\mbox{if }
i=j$;\vspace*{2pt}
\cr
-a_{i-1,i-j},&\quad $\mbox{if } 2\leq i\leq n, 1\leq j\leq
i-1$. }
\]
Hence,
%
\begin{equation}
\bolds{\Omega}^{-1}_{n}=\mathbf{T}_{n}'
\mathbf{D}^{-1}_{n}\mathbf {T}_{n}.\label{invautoCov}
\end{equation}
Because there are too many parameters in $\mathbf{T}_{n}$ and $\mathbf{D}_{n}$,
we are led to consider a banded Cholesky decomposition of
$\bolds{\Omega}^{-1}_{n}$,
%
\begin{equation}
\bolds{\Omega}^{-1}_{n}(k)=\mathbf{T}_{n}'(k)
\mathbf {D}^{-1}_{n}(k)\mathbf{T}_{n}(k),\label{invautoCov1}
\end{equation}
where $1\leq k\ll n$ is referred to as the banding parameter and
allowed to grow to infinity with $n$,
\[
\mathbf{D}_{n}(k)=\operatorname{diag}\bigl(\gamma_{0},
\sigma_{1}^{2},\ldots, \sigma _{k}^{2},
\ldots, \sigma_{k}^{2}\bigr),
\]
and
$\mathbf{T}_{n}(k)=(t_{ij}(k))_{1\leq i,j\leq n}$ with
\[
t_{ij}(k)= \cases{ 0,&\quad$\mbox{if } i < j \mbox{ or } \{ k+1 < i \leq n,
1 \leq j \leq i-k-1 \}$;\vspace*{2pt}
\cr
1,&\quad$\mbox{if } i=j$;\vspace*{2pt}
\cr
-a_{i-1,i-j},&\quad$\mbox{if } 2\leq i\leq k, 1\leq j\leq i-1$;\vspace *{2pt}
\cr
-a_{k,i-j},&\quad $\mbox{if }k+1 \leq i \leq n, i-k \leq j \leq i-1$. }
\]
We propose estimating $\bolds{\Omega}^{-1}_{n}$ using the sample
counterpart of
\eqref{invautoCov1},
%
\begin{equation}
\hat{\bolds{\Omega}}{}^{-1}_{n}(k):=\hat{
\mathbf{T}}_{n}'(k)\hat {\mathbf{D}}{}^{-1}_{n}(k)
\hat{\mathbf{T}}_{n}(k),\label{invautoCovest}
\end{equation}
where $\hat{\mathbf{T}}_{n}(k)$ and $\hat{\mathbf{D}}_{n}(k)$ are
obtained by plugging
in the least squares estimates of
the coefficients in $\mathbf{T}_{n}(k)$
and the corresponding residual variances in $\mathbf{D}_{n}(k)$; see
Section~\ref{section3} for more details.

Under \eqref{model}--\eqref{model-2}, this paper establishes
%
\begin{eqnarray}
\bigl\|\hat{\bolds{\Omega}}{}^{-1}_{n}(k)-\bolds{
\Omega}^{-1}_{n}\bigr\| _{2}=\mathrm{o}_{p}(1),
\label{keyresult}
\end{eqnarray}
with $k=K_{n} \to\infty$ satisfying \eqref{Kn}.
To appreciate the subtlety of \eqref{keyresult},
note that if $m$ independent realizations
$\mathbf{U}^{(1)}=(u^{(1)}_{1}, \ldots,u^{(1)}_{n})^{\prime}, \ldots,
\mathbf{U}^{(m)}=(u^{(m)}_{1}, \ldots,u^{(m)}_{n})^{\prime}$
of $\{u_{t}\}$
are available, Bickel and Levina \cite{r3}
introduced alternative estimates
$\check{\mathbf{T}}_{n, m}(k)$ and $\check{\mathbf{D}}_{n, m}(k)$
of $\mathbf{T}_{n}(k)$ and
$\mathbf{D}_{n}(k)$ through a multivariate analysis approach,
where $k<m<n$. More specifically, set $\tilde{\mathbf
{U}}_{j}=(u^{(1)}_{j}, \ldots, u^{(m)}_{j})^{\prime}$
and denote the regression coefficients of $\tilde{\mathbf{U}}_{j}$
on $\tilde{\mathbf{U}}_{j-1}, \ldots, \tilde{\mathbf{U}}_{\max{\{j-k,
1\}}}$
by $\check{\mathbf{a}}_{j}$.
Then $\check{\mathbf{T}}_{n, m}(k)$ and
$\check{\mathbf{D}}_{n, m}(k)$, respectively,
are obtained by replacing the coefficients in the $i$th row of
$\mathbf{T}_{n}(k)$ with $-\check{\mathbf{a}}_{i}$,
and $i$th diagonal element
of $\mathbf{D}_{n}(k)$
with the corresponding residual variance,
where $i=2,\ldots, n$.
Bickel and Levina \cite{r3} also showed that the resultant estimate
$\check{\bolds{\Omega}}{}^{-1}_{n, m}(k)=\check{\mathbf{T}}^{\prime}_{n, m}(k)
\check{\mathbf{D}}{}^{-1}_{n, m}(k)\check{\mathbf{T}}_{n, m}(k)$ of
$\bolds{\Omega}^{-1}_{n}$
has the property
%
\begin{eqnarray}
\bigl\|\check{\bolds{\Omega}}{}^{-1}_{n, m}(k)-\bolds{
\Omega }{}^{-1}_{n}\bigr\|_{2}=\mathrm{o}_{p}(1),
\label{keyresult2}
\end{eqnarray}
under $m \to\infty$, $m^{-1}\log n \to0$, $k=K_{n, m}\asymp(m/\log
n)^{\theta}$ for some $0<\theta<1/2$,
\eqref{lemma3.2-3},
and
%
\begin{equation}
\sup_{k\geq
1}\|\bolds{\Omega}_{k}
\|_{2}< \infty,\label{maxe}
\end{equation}
where $g(x)\asymp h(x)$ means that there exists a constant $0<C<\infty$
such that
$C\leq \liminf_{x\to\infty} h(x)/g(x)\leq\limsup_{x\to\infty}
h(x)/g(x) \leq C^{-1}$.
Since \eqref{maxe} fails to hold for long-memory processes like
\eqref{farima} and $m \to\infty$ is needed in \eqref{keyresult2},
the most distinctive feature of \eqref{keyresult} is that it holds
for one ($m=1$) realization, without imposing \eqref{maxe}. It is
also noteworthy that Cai, Ren and Zhou \cite{r41} have recently
established the
optimal rate of convergence for estimating the inverse of a Toeplitz
covariance matrix under the spectral norm. However, the covariance
matrix associated with (\ref{farima}) is still precluded by the
class of matrices considered in their paper, which needs to obey
assumptions like (\ref{maxe}) and (\ref{lemma3.2-3}).

The rest of the paper is organized as follows.
In Section~\ref{section2}, we analyze
the difference between $\bolds{\Omega}^{-1}_{n}(k)$
and $\bolds{\Omega}^{-1}_{n}$.
In particular, by deriving convergence rates of
$\|\mathbf{T}_{n}(k)-\mathbf{T}_{n}\|_{2}$ and $\|\mathbf
{D}_{n}(k)-\mathbf{D}_{n}\|_{2}$,
we obtain a convergence rate of $\|\bolds{\Omega
}^{-1}_{n}(k)-\bolds{\Omega}^{-1}_{n}\|_{2}$,
which plays an indispensable role in the proof of \eqref{keyresult}.
Section~\ref{section3} is devoted to proving \eqref{keyresult}.
By establishing a number of sharp bounds for the higher moments
of the quadratic forms in $u_{t}$, we obtain a convergence rate
of $\|\hat{\bolds{\Omega}}{}^{-1}_{n}(k)-\bolds{\Omega
}^{-1}_{n}(k)\|_{2}$,
which, in conjunction with the results in Section~\ref{section2},
leads to
a convergence rate of
$\|\hat{\bolds{\Omega}}{}^{-1}_{n}(k)-\bolds{\Omega}^{-1}_{n}\|_{2}$,
and hence
\eqref{keyresult}.
In Section~\ref{section4}, the results in Section~\ref{section3} are
extended to
regression models with long-memory
errors satisfying \eqref{model}--\eqref{model-2}.
Specifically, we show that
when the unobservable long-memory errors are replaced by
the corresponding least squares residuals,
our estimate of $\bolds{\Omega}^{-1}_{n}$
still has the same convergence rate,
without imposing \textit{any} assumptions on the design matrices.
Moreover, the estimated matrix
is applied to construct
an estimate of
the finite-past predictor coefficient vector
of the error process,
and an FGLSE of the regression coefficient vector.
Rates of convergence
of the latter two estimates
are also derived in a somewhat intricate way.
In Section~\ref{section5}, we present a Monte Carlo study of the
finite-sample performance of the proposed inverse matrix estimates.

\section{Bias analysis of banded Cholesky factors}\label{section2}
Our analysis of $\|\bolds{\Omega}^{-1}_{n}-\bolds{\Omega
}^{-1}_{n}(k)\|_{2}$
is reliant on the following two conditions on $a_{m,i}$'s defined in~\eqref{regcoe}.
\begin{longlist}[(ii)]
\item[(i)]
There exists $C_{1}>0$ such that for any $1\leq i\leq m<\infty$,
%
\begin{equation}
\biggl\llvert \frac{a_{m,i}}{a_{i}}\biggr\rrvert \leq C_{1} \biggl(
\frac{m}{m-i+1} \biggr)^{d}.\label{yuleCon1}
\end{equation}
\item[(ii)]
There exist $C_{2}>0$ and $0<\delta<1$ such that for any $1\leq i\leq
\delta m$ and $1\leq m <\infty$,
%
\begin{equation}
\biggl\llvert \frac{a_{m,i}}{a_{i}}-1\biggr\rrvert \leq C_{2}
\frac{i}{m}.\label{yuleCon2}
\end{equation}
\end{longlist}

Some comments on \eqref{yuleCon1} and \eqref{yuleCon2} are in order.
Note first that \eqref{yuleCon1} and \eqref{yuleCon2}, together with~\eqref{model-2}, immediately imply
that there exists $C>0$ such that for any $k=1,2,\ldots,$
%
\begin{equation}
\sum_{i=1}^{k}|a_{k,i}|\leq C
\label{resu1},
\end{equation}
which will be used frequently in the sequel.
Throughout the rest of the paper, $C$ denotes a generic positive
constant independent of
any unbounded index sets of positive integers.
These two conditions assert that the
finite-past predictor coefficients $a_{m,i}, i=1, \ldots, m $
approach to the corresponding infinite-past predictor coefficients
$a_{1}, a_{2}, \ldots$
in a nonuniform way.
More specifically, they require that
$a_{m,i}/a_{i}$ is very close to 1 when $i=\mathrm{o}(m)$, but has order of magnitude
$m^{(1-\theta)d}$ when $m-i\asymp m^{\theta}$ with $0 \leq\theta<1$.
This does not seem to be counterintuitive
because for a long-memory process, the finite order truncation
tends to create severer upward distortions in those
$a_{i}$'s with $i$ near the truncation lag $m+1$.
In fact, when $\{u_{t}\}$ is an $\mathrm{I}(d)$ process with $0<d<1/2$,
\eqref{yuleCon1} and \eqref{yuleCon2} follow directly from
the proof of Theorem~13.2.1 of Brockwell and Davis \cite{r4}.
In the following, we shall show that
\eqref{yuleCon1} and \eqref{yuleCon2} are satisfied
by model \eqref{farima}. To this end, we need an auxiliary lemma.

%
\begin{lemma}\label{lemma2.1}
Assume \eqref{model}, \eqref{model2},
%
\begin{eqnarray}
\psi_{i} \sim i^{-1+d}\mathcal{L}(i)\label{maCon}
\end{eqnarray}
and
%
\begin{eqnarray}
a_{i} \sim \frac{i^{-1-d}d\sin(\pi d)}{\pi\mathcal{L}(i)},\label{arCon}
\end{eqnarray}
where $g(x)\sim h(x)$ if $\lim_{x\to\infty}g(x)/h(x)=1$ and
$\mathcal{L}(x)$ is a positive slowly varying function,
namely, $\lim_{x \rightarrow\infty}\mathcal{L}(\lambda x)/\mathcal
{L}(x)=1$ for all $\lambda>0$.
Then for all large $m$,
\[
\max_{1\leq i\leq m}\biggl\llvert \frac{m(a_{m,i}-a_{i})}{\sum_{j=i\wedge
(m+1-i)}^{\infty}|a_{j}|}\biggr\rrvert \leq
C,
\]
where $u\wedge v= \min\{u, v\}$.
\end{lemma}

\begin{pf}
For $h,j\in\mathbf{N}\cup\{0\}$, we define
\[
d_{s}(h,j)= \cases{ \xi_{h+j}, &\quad$\mbox{if }s=1$;
\vspace*{2pt}
\cr
\displaystyle\sum_{v=0}^{\infty}
\xi_{h+j+v}d_{s-1}(h,v), &\quad$\mbox{if }s=2,3,\ldots,$ }
\]
where
$\xi_{t}=\sum_{v=0}^{\infty}\psi_{v}a_{v+t}$ for $t=0,1,\ldots.$
By Theorem~2.9 of Inoue and Kasahara \cite{r10}, we obtain
%
\begin{eqnarray}\label{lemma2.1-1}
&&m\bigl|(a_{m,i}-a_{i})\bigr|
\nonumber
\\[-8pt]
\\[-8pt]
\nonumber
&&\quad=m\Biggl\llvert \sum
_{s=1}^{\infty}\sum_{j=0}^{\infty}a_{i+j}d_{2s}(m+1,j)+
\sum_{s=1}^{\infty}\sum
_{j=0}^{\infty}a_{m+1-i+j}d_{2s-1}(m+1,j)
\Biggr\rrvert .
\end{eqnarray}
Let $\kappa>1$ satisfy
$0<\kappa\sin(\pi d)<1$. According to Proposition~3.2(i) of Inoue and
Kasahara \cite{r10},
there exists $N\in\mathbf{N}$ such that
%
\begin{equation}
0< d_{s}(h,j)\leq\frac{g_{s}(0)\{\kappa\sin(\pi d)\}^{s}}{h}, \qquad j\in\mathbf{N}\cup\{0\}, s\in
\mathbf{N}, h\geq N,\label{lemma2.1-2}
\end{equation}
where
\[
g_{s}(x)= \cases{ \displaystyle\frac{1}{\pi(1+x)}, &\quad$\mbox{if }s=1$;\vspace*{2pt}
\cr
\displaystyle\frac{1}{\pi^{2}}\int_{0}^{\infty}
\frac{\mathrm{d}v_{1}}{(v_{1}+1)(v_{1}+1+x)}, &\quad$\mbox{if }s=2$;\vspace*{2pt}
\cr
\displaystyle\frac{1}{\pi^{s}}\int
_{0}^{\infty}\cdots\int_{0}^{\infty}
\frac
{1}{v_{s-1}+1}\vspace*{2pt}\cr
\hspace*{56pt}{}\displaystyle\times \Biggl\{\prod_{j=1}^{s-2}
\frac
{1}{v_{j+1}+v_{j}+1} \Biggr\} \frac{1}{v_{1}+1+x}\,\mathrm{d}v_{s-1}\cdots\,
\mathrm{d}v_{1}, &\hspace*{-2pt}\quad$\mbox{if }s=3,4,\ldots.$ }
\]
Thus, for $m\geq N$ and $i=1,2,\ldots,m$,
%
\begin{eqnarray}\label{lemma2.1-3}
&&m\Biggl\llvert \sum_{s=1}^{\infty}\sum
_{j=0}^{\infty
}a_{i+j}d_{2s}(m+1,j)
\Biggr\rrvert\nonumber\\
&&\quad \leq m\sum_{j=0}^{\infty}|a_{i+j}|
\Biggl(\sum_{s=1}^{\infty}\frac
{g_{2s}(0)\{\kappa\sin(\pi d)\}^{2s}}{m+1}
\Biggr)
\\
&&\quad=\frac{m}{m+1}\sum_{j=i}^{\infty}|a_{j}|
\sum_{s=1}^{\infty}g_{2s}(0)\bigl\{
\kappa\sin(\pi d)\bigr\}^{2s}\leq C\sum_{j=i}^{\infty}|a_{j}|,
\nonumber
\end{eqnarray}
where the first inequality is by \eqref{lemma2.1-2} and the last one is by
Lemma~3.1(i) of Inoue and Kasahara~\cite{r10}. Similarly, \eqref
{lemma2.1-2} and Lemma~3.1(ii) of Inoue and Kasahara \cite{r10} imply
%
\begin{equation}
m\Biggl\llvert \sum_{s=1}^{\infty}\sum
_{j=0}^{\infty
}a_{m+1-i+j}d_{2s-1}(m+1,j)
\Biggr\rrvert \leq C\sum_{j=m+1-i}^{\infty
}|a_{j}|.\label{lemma2.1-4}
\end{equation}
Combining \eqref{lemma2.1-1}, \eqref{lemma2.1-3} and \eqref{lemma2.1-4}
yields the desired conclusion.
\end{pf}
%
%
%
\begin{rem}\label{remark2.1}
Theorem~3.3 of Inoue and Kasahara \cite{r10} shows that for any fixed
integer $i$,
%
\begin{equation}
\lim_{m\to\infty}m(a_{m,i}-a_{i})=d^{2}
\sum_{j=i}^{\infty}a_{j}.
\label{ik}
\end{equation}
Therefore, Lemma~\ref{lemma2.1} can be viewed as a \textit{uniform} extension
of \eqref{ik}.
\end{rem}
%
%
%
\begin{rem}\label{remark2.2}
Note that \eqref{model-1} and \eqref{model-2}
are fulfilled by \eqref{maCon} and \eqref{arCon}
if
%
\begin{equation}
0<\inf_{i}\mathcal{L}(i) \leq\sup_{i}
\mathcal{L}(i)<\infty. \label{nicefun}
\end{equation}
\end{rem}
%
By making use of Lemma~\ref{lemma2.1}, the next theorem shows that
\eqref{yuleCon1} and \eqref{yuleCon2} are met by
\eqref{maCon} and \eqref{arCon} with $\mathcal{L}(i)$
obeying \eqref{nicefun}.
Since the coefficients of the MA and AR representations of \eqref{farima}
take the form of \eqref{maCon} and \eqref{arCon}, respectively, for which
$\mathcal{L}(i)$ is a constant function (see Corollary~3.1 of Kokoszka
and Taqqu \cite{r32}
and Example~2.6 of Inoue and Kasahara \cite{r10}),
this theorem guarantees that
\eqref{farima} satisfies
\eqref{yuleCon1} and \eqref{yuleCon2}, confirming the flexibility of
these two conditions.
%
\begin{thmm}\label{theorem2.1}
Under the same assumptions as in Lemma~\ref{lemma2.1} with $\mathcal
{L}(i)$ satisfying \eqref{nicefun},
we have~\eqref{yuleCon1} and \eqref{yuleCon2}.
\end{thmm}
\begin{pf}
It suffices to show that \eqref{yuleCon1} and \eqref{yuleCon2}
hold for all sufficiently large $m$.
By Lemma~\ref{lemma2.1} and \eqref{nicefun}, it follows that for
all $1\leq i \leq m$ and all large $m$,
\[
\bigl|m(a_{m,i}-a_{i})\bigr|\leq C\max\bigl\{i^{-d},(m-i+1)^{-d}
\bigr\},
\]
yielding
\begin{eqnarray*}
\biggl\llvert \frac{a_{m,i}}{a_{i}}\biggr\rrvert =\biggl\llvert \frac{m(a_{m,i}-a_{i})}{ma_{i}}+1
\biggr\rrvert \leq C \frac{\max\{i^{-d},(m-i+1)^{-d}\}}{m i^{-1-d}}+1 \leq C \biggl(\frac{m}{m-i+1}
\biggr)^{d}.
\end{eqnarray*}
Therefore, \eqref{yuleCon1} follows.
Similarly, for all $1\leq i\leq\delta m$ with $0<\delta<1$ and all
large $m$,
\begin{eqnarray*}
\biggl\llvert \frac{a_{m,i}-a_{i}}{a_{i}}\biggr\rrvert &\leq& C \frac{\max\{
i^{-d},(m-i+1)^{-d}\}}{m i^{-1-d}} \leq C
\frac{i}{m},
\end{eqnarray*}
which leads to \eqref{yuleCon2}. Thus the proof is complete.
\end{pf}
%
Throughout the rest of this paper,
let $K_n$ denote a sequence of numbers satisfying
$K_n \rightarrow\infty$ and $K_n/n \rightarrow0$ as $n \rightarrow
\infty$. We are now ready to
provide upper bounds for
$\|\mathbf{T}_{n}-\mathbf{T}_{n}(K_{n})\|_{2}$
and $\|\mathbf{D}_{n}^{-1}-\mathbf{D}_{n}(K_{n})^{-1}\|_{2}$ in
Propositions \ref{proposition2.1}
and \ref{proposition2.2},
which in turn lead to a rate of convergence of
$\|\bolds{\Omega}^{-1}_{n}-\bolds{\Omega}^{-1}_{n}(K_{n})\|
_{2}$ in Theorem~\ref{theorem2.2}.
Before proceeding, we need two technical lemmas.

%
\begin{lemma}\label{lemma2.2}
Assume \eqref{model-2}, \eqref{yuleCon1} and \eqref{yuleCon2}.
Then
\begin{longlist}[(iii)]
\item[(i)]
$\sum_{j=K_{n}+1}^{k}|a_{k,j}|\leq C K_{n}^{-d}$ for any ${K_{n}+1\leq
k\leq n-1}$.
\item[(ii)]
$\sum_{j=1}^{K_{n}}|a_{k,j}-a_{K_{n},j}|\leq CK_{n}^{-d}$ for any
${K_{n}+1\leq k\leq n-1}$.
\item[(iii)]
$\sum_{j=\max(1,K_{n}+1-k)}^{K_{n}}|a_{j+k,j}-a_{K_{n},j}|\leq CK_{n}^{-d}$
for any ${1\leq k\leq n-K_{n}-1}$.
\item[(iv)]
$\sum_{j=1}^{n-k-1}|a_{j+k,j}-a_{K_{n},j}|\leq CK_{n}^{-d}$ for any
${n-K_{n}\leq k\leq n-2}$.
\end{longlist}
\end{lemma}

\begin{pf}
The proof is straightforward, and thus omitted.
\end{pf}
%
%
\begin{lemma}\label{lemma2.3}
Assume \eqref{model}--\eqref{model-2}.
Then for any $k \geq1$, $\sigma^{2}_{k}-\sigma^{2} \leq Ck^{-1}$,
where $\sigma^{2}_{k}$ is defined in~\eqref{sigmak}.
\end{lemma}
\begin{pf}
In view of \eqref{model2} and \eqref{regcoe}, it follows that for any
$k \geq1$,
$\sigma^{2}_{k}-\sigma^{2} \leq E(\sum_{j=k+1}^{\infty}a_{j}\times\break u_{t-j})^{2}$.
In addition, by \eqref{acf} (which is ensured by \eqref{model} and
\eqref{model-1}),
\eqref{model-2}, and Theorem~2.1 of Ing and Wei~\cite{r11}, one has for
any $k \geq1$ and $m \geq k+1$,
$E(\sum_{j=k+1}^{m}a_{j}u_{t-j})^{2} \leq C k^{-1}$, which, together
with the previous inequality,
gives the desired conclusion.
\end{pf}

%
\begin{prop}\label{proposition2.1}
Under the same assumptions as in Lemma~\ref{lemma2.2},

\begin{longlist}[(ii)]
\item[(i)]
$\|\mathbf{T}_{n}-\mathbf{T}_{n}(K_{n})\|_{2}=\mathrm{O}((K_{n}^{-d}\log
n)^{1/2})$.
\item[(ii)]
$\|\mathbf{T}_{n}(K_{n})\|_{2}=\mathrm{O}((\log K_{n})^{1/2})$.
\end{longlist}
\end{prop}
\begin{pf}
Let $\|\mathbf{B}\|_{k}=\max_{\|\mathbf{z}\|_{k}=1}\|\mathbf{Bz}\|_{k}$
denote the $k$-norm of an $h\times h$ matrix $\mathbf{B}$, where
$\|\mathbf{z}\|_{k}=(\sum_{i=1}^{h} |z_i|^k)^{1/k}$ is the $k$-norm of
the vector $\mathbf{z}=(z_1, \ldots, z_h)'$.
Then, by Lemma~\ref{lemma2.2}(i) and (ii),
\[
\bigl\|\mathbf{T}_{n}-\mathbf{T}_{n}(K_{n})
\bigr\|_{\infty} =\max_{K_{n}+1\leq i\leq n-1}\sum_{j=K_{n}+1}^{i}|a_{i,j}|+
\sum_{j=1}^{K_{n}}|a_{K_{n},j}-a_{i,j}|=\mathrm{O}
\bigl(K_{n}^{-d}\bigr).
\]
Moreover, $\|\mathbf{T}_{n}-\mathbf{T}_{n}(K_{n})\|_{1}$ is the maximum of
\[
\max_{0\leq k\leq
n-K_{n}-1} \Biggl\{\sum_{i=0}^{n-K_{n}-k-2}
|a_{K_{n}+1+i+k,K_{n}+1+i}|+\sum_{j=\max
(1,K_{n}+1-k)}^{K_{n}}|a_{j+k,j}-a_{K_{n},j}|
\Biggr\}\label{proposition2.1-1}
\]
and
$\max_{n-K_{n}\leq k\leq n-2}\sum_{j=1}^{n-k-1}|a_{j+k,j}-a_{K_{n},j}|$.
By \eqref{yuleCon1} and Lemma~\ref{lemma2.2}(iii) and (iv),
\begin{eqnarray*}
\max_{0\leq k\leq
n-K_{n}-1}\sum_{i=0}^{n-K_{n}-k-2}|a_{K_{n}+1+i+k,K_{n}+1+i}|&=&\mathrm{O}(
\log n),
\\
\max_{0\leq k\leq
n-K_{n}-1}\sum_{j=\max
(1,K_{n}+1-k)}^{K_{n}}|a_{j+k,j}-a_{K_{n},j}|&=&\mathrm{O}
\bigl(K^{-d}_{n}\bigr)
\end{eqnarray*}
and
\[
\max_{n-K_{n}\leq k\leq n-2}\sum_{j=1}^{n-k-1}|a_{j+k,j}-a_{K_{n},j}|=\mathrm{O}
\bigl(K^{-d}_{n}\bigr).
\]
Hence,
$\|\mathbf{T}_{n}-\mathbf{T}_{n}(K_{n})\|_{1}=\mathrm{O}(\log n)$. The proof of
(i) is completed by
\[
\bigl\|\mathbf{T}_{n}-\mathbf{T}_{n}(K_{n})
\bigr\|_{2} \leq\bigl(\bigl\|\mathbf{T}_{n}-\mathbf{T}_{n}(K_{n})
\bigr\|_{1}\bigl\|\mathbf {T}_{n}-\mathbf{T}_{n}(K_{n})
\bigr\|_{\infty}\bigr)^{1/2}.
\]
Similarly, it can be shown that
$\|\mathbf{T}_{n}(K_{n})\|_{\infty} =\mathrm{O}(1)$
and $\|\mathbf{T}_{n}(K_{n})\|_{1}=\mathrm{O}(\log K_{n})$,
yielding (ii).
\end{pf}
%
%
\begin{prop}\label{proposition2.2}
Under the same assumptions as in Lemma~\ref{lemma2.3},

\begin{longlist}[(ii)]
\item[(i)]
$\|\mathbf{D}_{n}^{-1}-\mathbf{D}_{n}^{-1}(K_{n})\|_{2}=\mathrm{O}(K_{n}^{-1})$.
\item[(ii)]
$\|\mathbf{D}_{n}^{-1}(K_{n})\|_{2}=\mathrm{O}(1)$.
\end{longlist}
\end{prop}
\begin{pf}
Equation (i) is an immediate consequence of Lemma~\ref{lemma2.3}.
Equation (ii) follows from
$\|\mathbf{D}_{n}^{-1}(K_{n})\|_{2}=\max_{0
\leq k \leq K_{n}}\sigma_{k}^{-2} \leq\sigma^{-2}$,
where $\sigma^{2}_{0}=\gamma_0$.
\end{pf}
%
%
\begin{thmm}\label{theorem2.2}
Assume \eqref{model}--\eqref{model-2}, \eqref{yuleCon1} and \eqref{yuleCon2}.
Suppose
%
\begin{equation}
\frac{\log n\log K_{n}}{K_{n}^{d}}=\mathrm{o}(1).\label{rate0}
\end{equation}
Then
%
\begin{equation}
\bigl\|\bolds{\Omega}_{n}^{-1}-\bolds{
\Omega}_{n}^{-1}(K_{n})\bigr\| _{2}=\mathrm{O}\bigl(
\log n K_{n}^{-d}\log K_{n}\bigr)^{1/2}=\mathrm{o}(1).
\label{rate1}
\end{equation}
Moreover, if \eqref{model3} is assumed,
%
\begin{equation}
\bigl\|\bolds{\Omega}_{n}^{-1}(K_{n})
\bigr\|_{2}=\mathrm{O}(1). \label{rate2}
\end{equation}
\end{thmm}
\begin{pf}
Equation \eqref{rate1} follows directly from
Propositions
\ref{proposition2.1} and \ref{proposition2.2},
\begin{eqnarray*}
\bigl\|\bolds{\Omega}_{n}^{-1}-\bolds{
\Omega}_{n}^{-1}(K_{n})\bigr\|_{2} &\leq&\bigl\|
\mathbf{T}_{n}-\mathbf{T}_{n}(K_{n})
\bigr\|_{2}\bigl\|\mathbf{D}_{n}^{-1}\bigr\| _{2}\bigl(
\bigl\|\mathbf{T}_{n}-\mathbf{T}_{n}(K_{n})
\bigr\|_{2}+\bigl\|\mathbf {T}_{n}(K_{n})\bigr\|_{2}
\bigr)
\\
&&{}+\bigl\|\mathbf{T}_{n}(K_{n})\bigr\|_{2}\bigl\|
\mathbf{D}_{n}^{-1}-\mathbf {D}_{n}^{-1}(K_{n})
\bigr\|_{2}\bigl(\bigl\|\mathbf{T}_{n}-\mathbf{T}_{n}(K_{n})
\bigr\| _{2}+\bigl\|\mathbf{T}_{n}(K_{n})\bigr\|_{2}
\bigr)
\\
&&{}+\bigl\|\mathbf{T}_{n}(K_{n})\bigr\|_{2}\bigl\|
\mathbf{D}_{n}^{-1}(K_{n})\bigr\|_{2}\bigl\|
\mathbf{T}_{n}-\mathbf{T}_{n}(K_{n})
\bigr\|_{2},
\end{eqnarray*}
and \eqref{rate0}. 
Equations \eqref{rate1} and \eqref{lemma3.2-3} (which is ensured by
\eqref{model3}) further
lead to \eqref{rate2}.
\end{pf}
%
\section{Main results}\label{section3}
In the sequel, the following assumptions on the innovation process $\{
w_{t}\}$ of \eqref{model}
are frequently used:
\begin{longlist}[(M1)]
\item[(M1)] $\{w_{t}, \mathcal{F}_{t}\}$
is a martingale difference sequence, where
$\mathcal{F}_{t}$ is an
increasing sequence of $\sigma$-field generated by
$w_{s}, s \leq t$.
\item[(M2)]
$E(w_{t}^{2}|\mathcal{F}_{t-1})=\sigma^{2}$ a.s.
\item[(M3)]
For some $q\geq1$, there is a constant $C_{q}>0$ such that
\[
\sup_{-\infty<
t<\infty}E\bigl(|w_{t}|^{4q}|
\mathcal{F}_{t-1}\bigr)\leq C_{q}\qquad \mbox{a.s.}
\]
\end{longlist}
As mentioned in Section~\ref{section1}, $\hat{\mathbf{T}}_{n}(K_{n})$
is obtained by replacing
$\mathbf{a}(k)=(a_{k,1}, \ldots, a_{k,k})^{\prime}$
in $\mathbf{T}_{n}(K_{n})$ with the corresponding
the least squares estimates $\hat{\mathbf{a}}(k)=(\hat{a}_{k,1},\ldots
,\hat{a}_{k,k})^{\prime}$, where $k=1, \ldots, K_{n}$ and
\[
(\hat{a}_{k,1},\ldots,\hat{a}_{k,k})=\mathop{\operatorname{arg\,min}}_{(\alpha
_{1},\ldots,\alpha_{k})\in
R^{k}}\sum_{t=k+1}^{n}(u_{t}-
\alpha_{1}u_{t-1}-\alpha_{2}u_{t-2}-
\cdots-\alpha_{k}u_{t-k})^2.
\]
Similarly, $\hat{\mathbf{D}}_{n}(K_{n})$ is obtained by replacing
$\sigma_{k}^{2}$
in $\mathbf{D}_{n}(K_{n})$ with $\hat{\sigma}_{k}^{2}$, where $k=0,
\ldots, K_{n}$ and
\begin{eqnarray*}
\hat{\sigma}_{0}^{2}&=&(n-1)^{-1}\sum
_{t=1}^{n}(u_{t}-\bar{u})^{2},\qquad
\bar{u}=n^{-1}\sum_{t=1}^{n}u_{i},
\\
\hat{\sigma}_{k}^{2}&=&(n-k)^{-1}\sum
_{t=k+1}^{n} \Biggl(u_{t}-\sum
_{j=1}^{k} \hat{a}_{k,j}u_{t-j}
\Biggr)^{2}.
\end{eqnarray*}
Recall $\hat{\bolds{\Omega}}{}^{-1}_{n}(K_{n})=\hat{\mathbf
{T}}_{n}'(K_{n})\hat{\mathbf{D}}{}^{-1}_{n}(K_{n})\hat{\mathbf{T}}_{n}(K_{n})$.
The objective of this section is to show that
$\|\hat{\bolds{\Omega}}{}^{-1}_{n}(K_{n})-\bolds{\Omega
}_{n}^{-1}\|_{2}=\mathrm{o}_{p}(1)$ in Theorem~\ref{theorem3.1}. To this end, we develop rates of convergence of
$\|\hat{\mathbf{T}}_{n}(K_{n})-\mathbf{T}_{n}(K_{n})\|_{2}$ and
$\|\hat{\mathbf{D}}{}^{-1}_{n}(K_{n})-\mathbf{D}_{n}^{-1}(K_{n})\|_{2}$
in Propositions \ref{proposition3.1} and
\ref{proposition3.2}, respectively,
whose proofs are heavily reliant on the following four lemmas,
Lemmas \ref{lemma3.1}--\ref{lemma3.4}.

%
\begin{lemma}\label{lemma3.1}
Assume \eqref{model}--\eqref{model-2} and \textup{(M1)--(M3)}.
Let $\mathbf{U}_{t}(k)=(u_{t}, u_{t-1},\ldots, u_{t-k+1})'$ and
$w_{k,t+1}=u_{t+1}-\mathbf{a}(k)'\mathbf{U}_{t}(k)$. Then for any
$1\leq k \leq n-1$,
%
\begin{equation}
E\Biggl\llVert \frac{1}{n-k}\sum_{t=k}^{n-1}
\mathbf {U}_{t}(k) (w_{k,t+1}-w_{t+1})\Biggr\rrVert
_{2}^{2q}\leq C \biggl(\frac{1}{n-k}
\biggr)^{q(1-2d)} \label{T1}
\end{equation}
and
%
\begin{equation}
E\Biggl\llVert \frac{1}{n-k}\sum_{t=k}^{n-1}
\mathbf{U}_{t}(k)w_{k,t+1}\Biggr\rrVert _{2}^{2q}
\leq C \biggl\{ \biggl(\frac{1}{n-k} \biggr)^{q(1-2d)}+ \biggl(
\frac{k}{n-k} \biggr)^{q} \biggr\} \label{T2}.
\end{equation}
Moreover, for $\theta> 1/q$,
%
\begin{equation}
\max_{1\leq k \leq K_{n}}\Biggl\llVert \frac{1}{n-k}\sum
_{t=k}^{n-1}\mathbf {U}_{t}(k)w_{k,t+1}
\Biggr\rrVert _{2}^{2}= \mathrm{O}_{p} \biggl(
\frac{K^{\theta}_{n}}{n^{1-2d}}+\frac{K_{n}^{1+\theta
}}{n} \biggr). \label{T3}
\end{equation}
\end{lemma}
\begin{pf}
By \eqref{acf}, Lemma~\ref{lemma2.3} and an argument similar to
that used in Lemma~3 of Ing and Wei~\cite{r9}, one has for any $1\leq k
\leq n-1$,
\begin{eqnarray*}\label{lemma3.1-10}
&&E\Biggl\llVert \frac{1}{n-k}\sum_{t=k}^{n-1}
\mathbf {U}_{t}(k) (w_{k,t+1}-w_{t+1})\Biggr\rrVert
_{2}^{2q}
\\
&&\quad\leq C (n-k)^{-q}k^{q} \Biggl\{\bigl(
\sigma_{k}^{2}-\sigma^{2}\bigr)\sum
_{i=-(n-k)+1}^{n-k-1}|\gamma_{i}| \Biggr
\}^{q} \leq C \biggl(\frac
{1}{n-k} \biggr)^{q(1-2d)},
\end{eqnarray*}
which gives \eqref{T1}. Equation \eqref{T2} follows
from \eqref{T1}
and for any $1\leq k \leq n-1$,
%
\begin{equation}
E\Biggl\llVert \frac{1}{n-k}\sum_{t=k}^{n-1}
\mathbf{U}_{t}(k)w_{t+1}\Biggr\rrVert _{2}^{2q}
\leq Ck^{q}(n-k)^{-q},\label{lemma3.1-2}
\end{equation}
whose proof is exactly same as that of Lemma~4 of Ing and Wei \cite{r9}.
To show \eqref{T3}, note that by \eqref{T2} and $K_{n}=\mathrm{o}(n)$,
\begin{eqnarray*}
&& E\max_{1\leq k \leq K_{n}}\Biggl\llVert \frac{1}{n-k}\sum
_{t=k}^{n-1}\mathbf {U}_{t}(k)w_{k,t+1}
\Biggr\rrVert _{2}^{2q}
\\
&&\quad\leq C \sum_{k=1}^{K_{n}}\bigl
\{n^{-q(1-2d)}+k^{q}n^{-q}\bigr\} \leq C\bigl\{
K_{n}n^{-q(1-2d)}+K_{n}^{q+1}n^{-q}
\bigr\}.
\end{eqnarray*}
This, together with $\theta>1/q$, gives the desired conclusion \eqref{T3}.
\end{pf}

%
\begin{rem}\label{remark3.0}
Lemma A.1 of Godet \cite{r6} establishes
an inequality closely related to
\eqref{T1}.
In particular, the inequality
yields
\[
E\Biggl\llVert \frac{1}{\sqrt{n-k}}\sum_{t=k}^{n-1}
\mathbf {U}_{t}(k) (w_{k,t+1}-w_{t+1})\Biggr\rrVert
_{2}^{2q}\leq C\bigl\{k(n-k)^{2d}\bigr
\}^{q}\bigl(\sigma^{2}_{k}-\sigma^{2}
\bigr)^{q}.
\]
This bound together with Lemma~\ref{lemma2.3}
also leads to \eqref{T1}.
\end{rem}

%
\begin{lemma}\label{lemma3.2}
Let
\[
\hat{\bolds{\Gamma}}_{k,n}=\frac{1}{n-k}\sum
_{t=k}^{n-1}\mathbf {U}_{t}(k)
\mathbf{U}_{t}(k)'.
\]
Assume \eqref{model}, \eqref{model-1}, (\ref{model3})
and \textup{(M1)--(M3)} with $q=1$.
Suppose
%
\begin{equation}
K_{n}= \cases{ \mathrm{o}\bigl(n^{1/2}\bigr),&\quad$\mbox{if }0<d<1/4$;
\vspace*{2pt}
\cr
\mathrm{o}\bigl((n/\log n)^{1/2}\bigr),&\quad$\mbox{if }d=1/4$;
\vspace*{2pt}
\cr
\mathrm{o}\bigl(n^{1-2d}\bigr),&\quad$\mbox{if }1/4<d<1/2$.
}\label{conditionP}
\end{equation}
Then
%
\begin{equation}
\bigl\|\hat{\bolds{\Gamma}}{}^{-1}_{K_{n},n}\bigr\|_{2}=\mathrm{O}_{p}(1).
\label{T4}
\end{equation}
\end{lemma}
\begin{pf}
By the first moment bound theorem of Findley and Wei \cite{r5},
\eqref{acf} and an argument similar to that used in Lemma~2 of
Ing and Wei \cite{r9}, it follows that
%
\begin{eqnarray}
E\llVert \hat{\bolds{\Gamma}}_{K_{n},n}-\bolds{\Omega
}_{K_{n}}\rrVert _{2}^{2} =\cases{ \mathrm{O}
\bigl(K_{n}^{2}(n-K_{n})^{-1}\bigr),&\quad$
\mbox{if } 0<d<1/4$;\vspace*{2pt}
\cr
\mathrm{O}\bigl(K_{n}^{2}(n-K_{n})^{-1}
\log(n-K_{n})\bigr),&\quad$\mbox{if } d=1/4$;\vspace *{2pt}
\cr
\mathrm{O}
\bigl(K_{n}^{2}(n-K_{n})^{-2+4d}\bigr),&\quad$
\mbox{if } 1/4<d<1/2$. } \label{T4.1}
\end{eqnarray}
Combining this, \eqref{conditionP} and \eqref{model3} leads to \eqref{T4}.
\end{pf}

%
\begin{lemma}\label{lemma3.3}
Under the same assumptions as in Theorem~\ref{theorem2.2}, one has
for any $k \geq1$ and $m=0, \pm1, \pm2, \ldots,$
$\gamma_{\tau_{k}}(m)=C(|m|+1)^{-1+2d}$,
where with $\tau_{k, t}=u_{t+1}-w_{t+1}-\mathbf{a}'(k)\mathbf
{U}_{t}(k)=w_{k, t+1}-w_{t+1}$,
$\gamma_{\tau_{k}}(m)=E(\tau_{k, 1}\tau_{k, m+1})$.
\end{lemma}
\begin{pf}
This result follows by a tedious but direct calculation.
The details are omitted.
\end{pf}
%
%
\begin{lemma}\label{lemma3.4}
Assume that (\ref{yuleCon1}), (\ref{yuleCon2}), and the assumptions of
Lemma~\ref{lemma3.1} hold.
Then, for any $1\leq k \leq n-1$,
%
\begin{eqnarray}
E\Biggl\llvert \frac{1}{n-k}\sum_{t=k}^{n-1}w_{k,t+1}^{2}-
\sigma_{k}^{2}\Biggr\rrvert ^{2q}\leq \cases{
C(n-k)^{-q},&\quad$\mbox{if } 0<d<1/4$,\vspace*{2pt}
\cr
C\bigl(\bigl
\{(n-k)^{-1}\log(n-k)\bigr\}^{q}\bigr),&\quad$\mbox{if } d=1/4$,
\vspace*{2pt}
\cr
C(n-k)^{-2q+4qd},&\quad$\mbox{if } 1/4<d<1/2$. \label{T5} }
\end{eqnarray}
Moreover, for $\theta>1/(2q)$,
%
\begin{eqnarray}
\max_{1\leq k \leq K_{n}}\Biggl\llvert \frac{1}{n-k}\sum
_{t=k}^{n-1}w_{k,t+1}^{2}-
\sigma_{k}^{2}\Biggr\rrvert = \cases{ \mathrm{O}_{p}
\bigl(K^{\theta}_{n}n^{-1/2} \bigr),&\quad$\mbox{if }
0<d<1/4$,\vspace *{2pt}
\cr
\mathrm{O}_{p} \bigl(K^{\theta}_{n}
(\log n)^{1/2}n^{-1/2} \bigr),&\quad$\mbox{if } d=1/4$,\vspace*{2pt}
\cr
\mathrm{O}_{p} \bigl(K^{\theta}_{n}n^{-1+2d}
\bigr),&\quad$\mbox{if } 1/4<d<1/2$. \label{T6} }
\end{eqnarray}
\end{lemma}
\begin{pf}
To show \eqref{T5}, note first that
%
\begin{eqnarray}
E\Biggl\llvert \frac{1}{n-k}\sum_{t=k}^{n-1}w_{k,t+1}^{2}-
\sigma_{k}^{2}\Biggr\rrvert ^{2q} \leq C
\bigl(E\bigl|(A1)\bigr|^{2q}+E\bigl|(A2)\bigr|^{2q}+E\bigl|(A3)\bigr|^{2q}\bigr),
\label{T6.1}
\end{eqnarray}
where
$E|(A1)|^{2q}=E|\frac{1}{n-k}\sum_{t=k}^{n-1}w_{t+1}^{2}-\sigma^{2}|^{2q}$,
$E|(A2)|^{2q}=E|\frac{2}{n-k}\sum_{t=k}^{n-1}w_{t+1}\tau_{k,t}|^{2q}$,
and $E|(A3)|^{2q}=E|\frac{1}{n-k}\sum_{t=k}^{n-1}\tau_{k,t}^{2}-(\sigma
^{2}_{k}-\sigma^{2})|^{2q}$.
It is clear that for any $1\leq k \leq n-1$,
%
\begin{equation}
E\bigl|(A1)\bigr|^{2q}\leq C(n-k)^{-q}. \label{T6.2}
\end{equation}
In addition, the first moment bound theorem of Findley and Wei \cite
{r5} implies that for any $1\leq k \leq n-1$,
\begin{eqnarray*}
E\bigl|(A2)\bigr|^{2q}&\leq& C\bigl((n-k)^{-1}\gamma_{\tau_{k}}(0)
\bigr)^{q},
\\
E\bigl|(A3)\bigr|^{2q}&\leq& C \Biggl\{(n-k)^{-1}\sum
_{j=0}^{n-k-1}\gamma_{\tau_{k}}^{2}(j)
\Biggr\}^{q},
\end{eqnarray*}
which, together with Lemmas \ref{lemma2.3} and \ref{lemma3.3}, \eqref
{T6.1} and \eqref{T6.2}, yield \eqref{T5}.
Equation \eqref{T6} follows immediately from \eqref{T5}
and an argument similar to that used to prove \eqref{T3}. The details
are omitted.
\end{pf}

We are now ready to establish rates of convergence of $\|\hat{\mathbf
{T}}_{n}(K_{n})-\mathbf{T}_{n}(K_{n})\|_{2}$
and $\|\hat{\mathbf{D}}{}^{-1}_{n}(K_{n})-\mathbf{D}_{n}^{-1}(K_{n})\|_{2}$.

%
\begin{prop}\label{proposition3.1}
Assume \eqref{model}--\eqref{model-2}, \eqref{model3} and \textup{(M1)--(M3)}.
Suppose (\ref{conditionP}).
Then for any $\theta>1/q$,
%
\begin{equation}
\bigl\|\hat{\mathbf{T}}_{n}(K_{n})-\mathbf{T}_{n}(K_{n})
\bigr\|^{2}_{2}= \mathrm{O}_{p} \biggl(\frac{K_{n}^{1+\theta}}{n^{1-2d}}+
\frac{K^{2+\theta}_{n}}{n} \biggr). \label{T7}
\end{equation}
\end{prop}
\begin{pf}
Let
$\mathbf{S}_{n}=(s_{ij})_{1\leq i, j\leq n}=\hat{\mathbf
{T}}_{n}(K_{n})-\mathbf{T}_{n}(K_{n})$.
Then
\[
\max_{1\leq i \leq n} \sum_{t=1}^{n}s_{it}^{2}
\leq\max_{1\leq k \leq
K_{n}}\bigl\|\hat{\mathbf{a}}(k)-\mathbf{a}(k)
\bigr\|_{2}^{2},
\]
and for each $1\leq j \leq n$,
$\sharp B_j \leq2K_{n}-1$, where $B_j=\{i \dvt \sum_{t=1}^n s_{it}s_{jt}
\neq0\}$.
These and some algebraic manipulations yield
\begin{eqnarray*}
\bigl\|\mathbf{S}_{n}\mathbf{S}_{n}'
\bigr\|_{1}&=&\max_{1\leq j \leq n}\sum
_{i
\in B_j}\Biggl\llvert \sum_{t=1}^{n}s_{it}s_{jt}
\Biggr\rrvert
\\
&\leq& \max_{1\leq j \leq n}\sum_{i \in B_j}
\Biggl(\sum_{t=1}^{n}s_{it}^{2}
\Biggr)^{1/2} \Biggl(\sum_{h=1}^{n}s_{jh}^{2}
\Biggr)^{1/2} \leq CK_{n}\max_{1\leq k \leq
K_{n}}\bigl\|\hat{
\mathbf{a}}(k)-\mathbf{a}(k)\bigr\|_{2}^{2}
\\
&\leq& CK_{n}\bigl\|\hat{\bolds{\Gamma}}{}^{-1}_{K_{n},n}
\bigr\|_{2}^{2} \max_{1\leq k \leq
K_{n}} \Biggl\llVert
\frac{1}{n-k}\sum_{t=k}^{n-1}
\mathbf{U}_{t}(k)w_{k,t+1}\Biggr\rrVert _{2}^{2}.
\end{eqnarray*}
Now, the desired conclusion \eqref{T7} follows from
\eqref{T3}, \eqref{T4} and
$\|\mathbf{S}_{n}\|^{2}_{2} \leq\|\mathbf{S}_{n}\mathbf{S}_{n}'\|_{1}$.
\end{pf}

%
\begin{prop}\label{proposition3.2}
Assume (\ref{yuleCon1}), (\ref{yuleCon2}), and the same assumptions as
in Proposition~\ref{proposition3.1}.
Suppose~(\ref{conditionP}).
Then for any $\theta>1/q$,
%
\begin{equation}
\bigl\|\hat{\mathbf{D}}{}^{-1}_{n}(K_{n})-\mathbf
{D}_{n}^{-1}(K_{n})\bigr\|_{2}= \cases{
\mathrm{O}_{p}\bigl(n^{-1/2}K_{n}^{\theta}\bigr),&\quad$
\mbox{if } 0<d<1/4$,\vspace*{2pt}
\cr
\mathrm{O}_{p}\bigl((\log
n/n)^{1/2}K_{n}^{\theta}\bigr),&\quad$\mbox{if } d=1/4$,
\vspace *{2pt}
\cr
\mathrm{O}_{p}\bigl(n^{-1+2d}K_{n}^{\theta}
\bigr),&\quad $\mbox{if } 1/4<d<1/2$. } \label{T8}
\end{equation}
\end{prop}
\begin{pf}
Note first that
%
\begin{equation}
\bigl\llVert \hat{\mathbf{D}}_{n}(K_{n})-
\mathbf{D}_{n}(K_{n})\bigr\rrVert _{2}= \max
_{0\leq k\leq
K_{n}}\bigl|\hat{\sigma}_{k}^{2}-
\sigma_{k}^{2}\bigr|, \label{T9.1}
\end{equation}
recalling $\sigma^{2}_{0}=\gamma_{0}$.
By \eqref{acf} and an argument similar to that used to prove \eqref
{T4.1}, it holds that
%
\begin{equation}
E\bigl(\hat{\sigma}_{0}^{2}-\sigma^{2}_{0}
\bigr)^{2}= \cases{ \mathrm{O}\bigl(n^{-1}\bigr),&\quad$\mbox{if } 0<d<1/4$,
\vspace*{2pt}
\cr
\mathrm{O}(\log n/n),&\quad$\mbox{if } d=1/4$,\vspace*{2pt}
\cr
\mathrm{O}
\bigl(n^{-2+4d}\bigr),&\quad$\mbox{if } 1/4<d<1/2$. }\label{proposition3.2-1}
\end{equation}
Straightforward calculations show
\begin{eqnarray*}
\max_{1 \leq k \leq K_{n}}\bigl|\hat{\sigma}_{k}^{2}-
\sigma_{k}^{2}\bigr| &\leq& \max_{1 \leq k \leq K_{n}}\Biggl
\llvert \frac{1}{n-k}\sum_{t=k}^{n-1}w_{k,t+1}^{2}-
\sigma_{k}^{2}\Biggr\rrvert
\\
&&{}+ C\bigl\|\hat{\bolds{\Gamma}}{}^{-1}_{K_{n},n}\bigr\|_{2}
\max_{1 \leq k \leq K_{n}}\Biggl\llVert \frac{1}{n-k}\sum
_{t=k}^{n-1}\mathbf {U}_{t}(k)w_{k,t+1}
\Biggr\rrVert ^{2}_{2},
\end{eqnarray*}
which, in conjunction with \eqref{proposition3.2-1}, \eqref{T9.1},
\eqref{T3}, \eqref{T4} and \eqref{T6},
results in \eqref{T8}.
\end{pf}

The main results of this section is given as follows.
%
\begin{thmm}\label{theorem3.1}
Assume the same assumptions as in Proposition~\ref{proposition3.2}.
Suppose
%
\begin{equation}
\frac{\log n \log K_{n}}{K^{d}_{n}} + \frac{K_{n}^{1+\theta} \log K_{n}}{n^{1-2d}}+ \frac{K_{n}^{2+\theta} \log K_{n}}{n} =\mathrm{o}(1),
\label{Kn}
\end{equation}
for some $\theta>1/q$.
Then
%
\begin{eqnarray}\label
{thmm3.1}
&&\bigl\|\hat{\bolds{\Omega}}{}^{-1}_{n}(K_{n})-
\bolds {\Omega}_{n}^{-1}\bigr\|_{2}
\nonumber
\\
&&\quad=\mathrm{O}_{p} \biggl( \biggl(\frac{\log n \log
K_{n}}{K^{d}_{n}} \biggr)^{1/2} +
\biggl(\frac{K_{n}^{1+\theta} \log K_{n}}{n^{1-2d}}+ \frac{K_{n}^{2+\theta} \log K_{n}}{n} \biggr)^{1/2}
\biggr)
\\
&&\quad=\mathrm{o}_{p}(1)
\nonumber
\end{eqnarray}
and
%
\begin{equation}
\bigl\|\hat{\bolds{\Omega}}{}^{-1}_{n}(K_{n})
\bigr\|_{2}=\mathrm{O}_{p}(1). \label{thmm3.1a}
\end{equation}
\end{thmm}
\begin{pf}
By Propositions \ref{proposition3.1} and \ref{proposition3.2}, \eqref{model3},
(ii) of Proposition~\ref{proposition2.1}, (ii) of Proposition~\ref
{proposition2.2}
and
\begin{eqnarray*}
&&\bigl\|\hat{\bolds{\Omega}}{}^{-1}_{n}(K_{n})-
\bolds{\Omega }{}^{-1}_{n}(K_{n})
\bigr\|_{2}
\\
&&\quad\leq\bigl\|\hat{\mathbf{T}}_{n}(K_{n})-\mathbf{T}_{n}(K_{n})
\bigr\|_{2}\bigl\|\hat {\mathbf{D}}{}^{-1}_{n}(K_{n})
\bigr\|_{2}\bigl(\bigl\|\hat{\mathbf{T}}_{n}(K_{n})-\mathbf
{T}_{n}(K_{n})\bigr\|_{2}+\bigl\|\mathbf{T}_{n}(K_{n})
\bigr\|_{2}\bigr)
\\
&&\qquad{}+\bigl\|\mathbf{T}_{n}(K_{n})\bigr\|_{2}\bigl\|\hat{\mathbf
{D}}{}^{-1}_{n}(K_{n})-\mathbf{D}_{n}^{-1}(K_{n})
\bigr\|_{2}\bigl(\bigl\|\hat{\mathbf {T}}_{n}(K_{n})-
\mathbf{T}_{n}(K_{n})\bigr\|_{2}+\bigl\|
\mathbf{T}_{n}(K_{n})\bigr\| _{2}\bigr)
\\
&&\qquad{}+\bigl\|\mathbf{T}_{n}(K_{n})\bigr\|_{2}\bigl\|
\mathbf{D}_{n}^{-1}(K_{n})\bigr\|_{2}\bigl\|\hat {
\mathbf{T}}_{n}(K_{n})-\mathbf{T}_{n}(K_{n})
\bigr\|_{2},
\end{eqnarray*}
one obtains
\begin{eqnarray*}
\bigl\|\hat{\bolds{\Omega}}{}^{-1}_{n}(K_{n})-
\bolds{\Omega }_{n}^{-1}(K_{n})
\bigr\|_{2}=\mathrm{O}_{p} \biggl( \biggl(\frac{K_{n}^{1+\theta} \log K_{n}}{n^{1-2d}}+
\frac{K_{n}^{2+\theta} \log K_{n}}{n} \biggr)^{1/2} \biggr).
\end{eqnarray*}
This, together with \eqref{Kn} and Theorem~\ref{theorem2.2}, leads to
the desired conclusions \eqref{thmm3.1}
and \eqref{thmm3.1a}.
\end{pf}

%
\begin{rem}\label{remark3.1}
It would be interesting
to compare Theorem~\ref{theorem3.1}
with the moment bounds for $\hat{\bolds{\Gamma}}{}^{-1}_{K_{n},n}$
given by
Godet \cite{r6}.
If $\{u_{t}\}$ is a Gaussian process satisfying
\eqref{model}--\eqref{model-2} and \eqref{model3},
then Theorem~2.1 of Godet \cite{r6} yields that for
%
\begin{eqnarray}
K_{n}&=&\mathrm{O}\bigl(n^{\lambda}\bigr)\qquad \mbox{with } 0<\lambda< \min
\{1/2, 1-2d\}, \label{Knn}
\\
E\bigl\|\hat{\bolds{\Gamma}}{}^{-1}_{K_{n},n}-\bolds {
\Omega}^{-1}_{K_{n}}\bigr\|_{2}&=& \cases{ \mathrm{O}
\bigl(n^{-1/2}K_{n}\bigr),&\quad$\mbox{if } 0<d<1/4$,\vspace*{2pt}
\cr
\mathrm{O}\bigl((\log n/n)^{1/2}K_{n}\bigr),&\quad$\mbox{if } d=1/4$,
\vspace*{2pt}
\cr
\mathrm{O}\bigl(n^{-1+2d}K_{n}\bigr),&\quad$\mbox{if }
1/4<d<1/2$.} \label{T88}
\end{eqnarray}
One major difference between
$\hat{\bolds{\Omega}}{}^{-1}_{n}(K_{n})$ and
$\hat{\bolds{\Gamma}}{}^{-1}_{K_{n},n}$ is that the former aims
at estimating the inverse autocovariance matrix of all $n$
observations, $\bolds{\Omega}^{-1}_{n}$, but the latter only
focuses on that of $K_{n}$ consecutive observations,
$\bolds{\Omega}^{-1}_{K_{n}}$, with $K_{n}\ll n$. While
\eqref{T88} plays an important role in analyzing the mean squared
prediction error of the least squares predictor of $u_{n+1}$ based
on the AR($K_{n}$) model, $\hat{\bolds{\Gamma}}{}^{-1}_{K_{n},n}$
cannot be used in situations where consistent estimates of
$\bolds{\Omega}_{n}^{-1}$ are indispensable. See Section~\ref{section4-2} for some examples. Moreover, the convergence rate
of $\hat{\bolds{\Omega}}{}^{-1}_{n}(K_{n})$ is determined by not
only the estimation error
$\|\hat{\bolds{\Omega}}{}^{-1}_{n}(K_{n})-\bolds{\Omega
}_{n}^{-1}(K_{n})\|_{2}$,
but also the approximation error
$\|\bolds{\Omega}_{n}^{-1}(K_{n})-\bolds{\Omega}_{n}^{-1}\|_{2}$.
This latter type of error, however, is irrelevant to the convergence
rate of $\hat{\bolds{\Gamma}}{}^{-1}_{K_{n},n}$.
Finally, we note that \eqref{T88} gives a stronger mode of
convergence than \eqref{thmm3.1}, but at the expense of more
stringent assumptions on moments and distributions.
\end{rem}

\section{Some extensions}\label{section4}
Consider a linear regression model with serially correlated errors,
%
\begin{equation}
y_{t} = \mathbf{x}'_{t} \bolds{\beta} + u_{t} = {\sum_{i=1}^{p}}
{x}_{t i} \beta_{i} + u_{t}, \label{regression model}
\end{equation}
where
$\bolds{\beta}$ is an unknown coefficient vector, $\mathbf
{x}_{t}$'s are $p$-dimensional nonrandom input vectors and $u_{t}$'s
are unobservable random disturbances satisfying the long-memory
conditions described previously.
Having observed $\mathbf{y}_{n} = ( y_{1}, \ldots, y_{n} )'$ and
$\check{\mathbf{x}}_{nj} = ( x_{1 j}, \ldots, x_{n j} )'$, $1 \leq j
\leq p$,
it is natural to
estimate $\mathbf{u}_{n} = ( u_{1}, \ldots, u_{n} )'$ via the least
squares residuals
\[
\tilde{\mathbf{u}}_{n} = ( \tilde{u}_{1}, \ldots,
\tilde{u}_{n} )' = ( \mathbf{I}_{n} -
M_{np} ) \mathbf{y}_{n} = ( \mathbf{I}_{n} -
M_{np} ) \mathbf{u}_{n},
\]
where $\mathbf{I}_{n}$ is the $n \times n$ identity matrix, and
$M_{np}$ is the orthogonal projection matrix of
$\overline{\operatorname{sp}} \{ \check{\mathbf{x}}_{n1}, \ldots,
\check{\mathbf{x}}_{np} \}$, the closed span of $\{
\check{\mathbf{x}}_{n1}, \ldots, \check{\mathbf{x}}_{np} \}$.
Note that $\tilde{\mathbf{u}}_{n}$
is also known as a detrended time
series, in particular when $\mathbf{x}_{t}$ represents the trend or
seasonal component of $y_{t}$.
Let $\{ \check{\mathbf{q}}_{ni}=( \mathsf{q}_{1 i}, \ldots, \mathsf
{q}_{n i} )', i =
1, \ldots, r \}$, $1 \leq r \leq p$, be an orthonormal basis of
$\overline{\operatorname{sp}} \{ \check{\mathbf{x}}_{n1}, \ldots,
\check{\mathbf{x}}_{np} \}$. It is well known that $M_{np}
=\sum_{i=1}^{r} \check{\mathbf{q}}_{ni} \check{\mathbf{q}}'_{ni}$, and hence
with ${v}_{i } = \check{\mathbf{q}}'_{ni} \mathbf{u}_{n}$,
%
\begin{equation}
\tilde{\mathbf{u}}_{n} = \mathbf{u}_{n} - {\sum
_{i=1}^{r}} {v}_{i} \check{
\mathbf{q}}_{ni}. \label{QV}
\end{equation}
In Section~\ref{section4-1}, we shall show that the inverse
autocovariance matrix, $\bolds{\Omega}_{n}^{-1}$, of
$\mathbf{u}_{n}$ can still be consistently estimated by the modified
Cholesky decomposition method proposed in Section~\ref{section3}
with $\mathbf{u}_{n}$ replaced by
$\tilde{\mathbf{u}}_{n}$, which is denoted by $\tilde{\bolds{\Omega
}}{}^{-1}_{n}(K_{n})$.
We also show that
$\tilde{\bolds{\Omega}}{}^{-1}_{n}(K_{n})$ and
$\hat{\bolds{\Omega}}{}^{-1}_{n}(K_{n})$ share the same rate of convergence.
Moreover, we propose an estimate of
$\mathbf{a}(n)=(a_{n,1}, \ldots, a_{n,n} )^{\prime}$, the $n$-dimensional
finite predictor coefficient vector of $\{u_{t}\}$, based on $\tilde
{\bolds{\Omega}}{}^{-1}_{n}(K_{n})$,
and derive its convergence rate.
These asymptotic results
are obtained without
imposing \textit{any} assumptions on the design matrix $\mathbf
{X}_{n}=(\check{\mathbf{x}}_{n1}, \ldots,
\check{\mathbf{x}}_{np})$.
On the other hand, we assume that
$\mathbf{X}_{n}$ has a full rank in Section~\ref{section4-2},
and propose an FGLSE of $\bolds{\beta}$ based on $\tilde
{\bolds{\Omega}}{}^{-1}_{n}(K_{n})$.
The rate of convergence of the proposed FGLSE is also established in Section~\ref{section4-2}.

\subsection{Consistent estimates of \texorpdfstring{$\bolds{\Omega}_{n}^{-1}$}{Omega{n}{-1}} and
$\mathbf{a}(n)$ based on \texorpdfstring{$\tilde{\mathbf{u}}_{n}$}{tilde{u}{n}}}\label{section4-1}
Define
\[
\tilde{\bolds{\Omega}}{}^{-1}_{n}(K_{n}):=
\tilde{\mathbf {T}}_{n}(K_{n})'\tilde{
\mathbf{D}}{}^{-1}_{n}(K_{n})\tilde{
\mathbf{T}}_{n}(K_{n})
\]
where
$\tilde{\mathbf{T}}_{n}(K_{n})$ and
$\tilde{\mathbf{D}}_{n}(K_{n})$ are $\hat{\mathbf{T}}_{n}(K_{n})$ and
$\hat{\mathbf{D}}_{n}(K_{n})$
with $\hat{a}_{ij}$ and $\hat{\sigma}_{i}^{2}$, respectively, replaced
by $\tilde{a}_{ij}$
and $\tilde{\sigma}_{i}^{2}$ defined as follows:
\begin{eqnarray*}
(\tilde{a}_{k,1},\ldots,\tilde{a}_{k,k})&=&\mathop{\operatorname{arg\,min}}_{(\alpha
_{1},\ldots,\alpha_{k})\in
R^{k}}\sum_{t=k+1}^{n}(
\tilde{u}_{t}- \alpha_{1}\tilde{u}_{t-1}-
\alpha_{2}\tilde{u}_{t-2}-\cdots-\alpha _{k}
\tilde{u}_{t-k})^2,
\\
\tilde{\sigma}_{0}^{2}&=&(n-1)^{-1}\sum
_{t=1}^{n}(\tilde{u}_{t}-\bar {
\tilde{u}})^{2}, \qquad\bar{\tilde{u}}=n^{-1}\sum
_{t=1}^{n}\tilde{u}_{i},
\\
\tilde{\sigma}_{k}^{2}&=&(n-k)^{-1}\sum
_{t=k+1}^{n} \Biggl(\tilde {u}_{t}-\sum
_{j=1}^{k} \tilde{a}_{k,j}
\tilde{u}_{t-j} \Biggr)^{2}.
\end{eqnarray*}
By establishing
probability bounds for
$\|\tilde{\mathbf{T}}_{n}(K_{n})-\mathbf{T}_{n}(K_{n})\|_{2}$ and
$\|\tilde{\mathbf{D}}{}^{-1}_{n}(K_{n})-\mathbf{D}{}^{-1}_{n}(K_{n})\|_{2}$
in Proposition~\ref{proposition4.1},
we obtain the convergence rate of
$\|\tilde{\bolds{\Omega}}{}^{-1}_{n}(K_{n})-\bolds{\Omega
}{}^{-1}_{n}\|_{2}$ in
Theorem~\ref{theorem4.1}.
According to (\ref{QV}), $\mathbf{u}_{n}$ and
$\tilde{\mathbf{u}}_{n}$ differ by the vector
$\sum_{i=1}^{r} {v}_{i} \check{\mathbf{q}}_{ni}$,
whose entries are weighted sums of $u_{1},u_{2},\ldots,u_{n}$ with weights
$\mathsf{q}_{t_{1},i}\mathsf{q}_{t_{2},j}$ for some $1\leq t_{1},
t_{2}\leq n$ and $1\leq i,j \leq r$.
To explore the contributions of
$\sum_{i=1}^{r} {v}_{i} \check{\mathbf{q}}_{ni}$ to
$\|\tilde{\mathbf{T}}_{n}(k)-\mathbf{T}_{n}(k)\|_{2}$ and
$\|\tilde{\mathbf{D}}{}^{-1}_{n}(k)-\mathbf{D}_{n}^{-1}(k)\|_{2}$,
we need moment bounds for
the linear combinations of
$u_{i}$'s
and $\tau_{k,i}$'s, which are introduced in the following lemma.

%
\begin{lemma}\label{lemma4.1}
Let $c_{1},\ldots,c_{m}$ be any real numbers. Under the same
assumptions as in Lemma~\ref{lemma3.1},
%
\begin{equation}
E \Biggl(\sum_{i=1}^{m}c_{i}u_{i}
\Biggr)^{4q} \leq C \Biggl(\sum_{i=1}^{m}c_{i}^{2}
\Biggr)^{2q}m^{4qd}. \label{4.11}
\end{equation}
Moreover, if (\ref{yuleCon1}) and (\ref{yuleCon2}) also hold true, then
%
\begin{equation}
E \Biggl(\sum_{i=1}^{m}c_{i}
\tau_{k,i} \Biggr)^{4q}\leq C \Biggl(\sum
_{i=1}^{m}c_{i}^{2}
\Biggr)^{2q}m^{4qd}. \label{4.12}
\end{equation}
\end{lemma}
\begin{pf}
By Lemma~2 of Wei \cite{r15}, we have $E(\sum_{i=1}^{m}c_{i}u_{i})^{4q}
\leq C\{E(\sum_{i=1}^{m}c_{i}u_{i})^{2}\}^{2q}$.
Theorem~2.1 of Ing and Wei \cite{r11} and Jensen's inequality further yield
$E(\sum_{i=1}^{m}c_{i}u_{i})^{2} \leq C (\sum_{i=1}^{m}|c_{i}|^{2/(1+2d)})^{1+2d}\leq
C m^{2d}(\sum_{i=1}^{m}c_{i}^{2})$. Hence, (\ref{4.11}) follows.
Equation (\ref{4.12}) is ensured by Lemma~\ref{lemma3.3}
and an argument similar to that used to prove (\ref{4.11}).
\end{pf}

Equipped with Lemma~\ref{lemma4.1}, we can prove another
auxiliary lemma, which plays a key role in establishing Proposition~\ref{proposition4.1}. First, some notation:
$\tilde{w}_{k,t+1}=\tilde{u}_{t+1}-\mathbf{a}(k)'\tilde{\mathbf{U}}_{t}(k)$,
$\tilde{\mathbf{U}}_{t}(k)=(\tilde{u}_{t},\tilde{u}_{t-1},\ldots,\tilde
{u}_{t-k+1})'$,
$\tilde{\bolds{\Gamma}}_{k,n}=\frac{1}{n-k}\sum_{t=k}^{n-1}\tilde
{\mathbf{U}}_{t}(k)\tilde{\mathbf{U}}_{t}(k)'$,
$\mathbf{q}_{t}=(\mathsf{q}_{t,1},\mathsf{q}_{t,2},\ldots,\mathsf{q}_{t,r})'$,
$\mathbf{Q}_{t}(k)=(\mathbf{q}_{t},\mathbf{q}_{t-1},\ldots,\mathbf{q}_{t-k+1})'$
and $\mathbf{V}_{n}=(v_{1}, \ldots, v_{r})^{\prime}$.

%
\begin{lemma}\label{lemma4.2}

\begin{longlist}[(iii)]
\item[(i)] Assume that the same assumptions as in Lemma~\ref{lemma3.4}
hold. Then for $K_{n}=\mathrm{o}(n)$ and $\theta>1/q$,
\begin{eqnarray*}
\max_{1\leq k\leq K_{n}}\Biggl\llVert \frac{1}{n-k}\sum
_{t=k}^{n-1}\tilde {\mathbf{U}}_{t}(k)
\tilde{w}_{k,t+1}\Biggr\rrVert _{2}^{2}
=\mathrm{O}_{p} \biggl(\frac{K^{\theta}_{n}}{n^{1-2d}}+\frac{K_{n}^{1+\theta
}}{n} \biggr).
\end{eqnarray*}
\item[(ii)] Assume that the same assumptions as in Lemma~\ref{lemma3.2}
hold. Then for $K_{n}$ satisfying \eqref{conditionP},
$\|\tilde{\bolds{\Gamma}}{}^{-1}_{K_{n},n}\|_{2}=\mathrm{O}_{p}(1)$.
\item[(iii)] Assume that the same assumptions as in Lemma~\ref
{lemma3.4} hold. Then for $K_{n}=\mathrm{o}(n)$ and $\theta>1/(2q)$,
\begin{eqnarray*}
\max_{1\leq k\leq K_{n}}\Biggl\llvert \frac{1}{n-k}\sum
_{t=k}^{n-1}\tilde {w}_{k,t+1}^{2}-
\sigma_{k}^{2}\Biggr\rrvert = \cases{ \mathrm{O}_{p}
\bigl(K^{\theta}_{n}n^{-1/2}\bigr),&\quad$\mbox{if } 0<d<1/4$;
\vspace*{2pt}
\cr
\mathrm{O}_{p}\bigl(K^{\theta}_{n} (\log
n)^{1/2}n^{-1/2}\bigr),&\quad$\mbox{if } d=1/4$;\vspace*{2pt}
\cr
\mathrm{O}_{p}\bigl(K^{\theta}_{n}n^{-1+2d}\bigr),&\quad$
\mbox{if } 1/4<d<1/2$. }
\end{eqnarray*}
\end{longlist}
\end{lemma}
\begin{pf}
We begin by proving (i).
Define
$(B1)=\|\frac{1}{n-k}\sum_{t=k}^{n-1}\tilde{\mathbf{U}}_{t}(k)(\tilde
{w}_{k,t+1}-w_{t+1})\|_{2}^{2q}$
and $(B2)=\|\frac{1}{n-k}\sum_{t=k}^{n-1}\tilde{\mathbf
{U}}_{t}(k)w_{t+1}\|_{2}^{2q}$.
Straightforward calculations yield
%
\begin{eqnarray}
\Biggl\llVert \frac{1}{n-k}\sum_{t=k}^{n-1}
\tilde{\mathbf{U}}_{t}(k)\tilde {w}_{k,t+1}\Biggr\rrVert
_{2}^{2q}\leq C\bigl\{(B1)+(B2)\bigr\}, \label{lemma4.2-1}
\end{eqnarray}
%
\begin{eqnarray}
E(B1) \leq Cn^{-2q}k^{q-1}\sum_{j=1}^{k}
\bigl\{ E\bigl|(B3)\bigr|^{2q}+E\bigl|(B4)\bigr|^{2q}+E\bigl|(B5)\bigr|^{2q}
+E\bigl|(B6)\bigr|^{2q}
\bigr\} \label{lemma4.2-1.5}
\end{eqnarray}
and
%
\begin{eqnarray}
E(B2) \leq C \Biggl\{E\Biggl\llVert \frac{1}{n-k}\sum
_{t=k}^{n-1}\mathbf {U}_{t}(k)w_{t+1}
\Biggr\rrVert _{2}^{2q}+E(B7) \Biggr\}, \label{lemma4.2-2}
\end{eqnarray}
where
\begin{eqnarray*}
(B3)&=&\sum_{t=k}^{n-1}u_{t+1-j}
\tau_{k,t}, \qquad (B4)=\mathbf{V}^{\prime}_{n}\sum
_{t=k}^{n-1}\mathbf{q}_{t+1-j}
\tau_{k,t},
\\
(B5)&=&\mathbf{V}^{\prime}_{n}\sum_{t=k}^{n-1}
\bigl(\mathbf{q}_{t+1}-\mathbf {Q}^{\prime}_{t}(k)
\mathbf{a}(k)\bigr)u_{t+1-j},
\\
(B6)&=&\mathbf{V}^{\prime}_{n} \Biggl\{\sum
_{t=k}^{n-1}\mathbf {q}_{t+1-j}\bigl(
\mathbf{q}_{t+1}'-\mathbf{a}(k)'
\mathbf{Q}_{t}(k)\bigr) \Biggr\} \mathbf{V}_{n},
\\
(B7)&=&\Biggl\llVert \Biggl\{\frac{1}{n-k}\sum_{t=k}^{n-1}
\mathbf {Q}_{t}(k)w_{t+1} \Biggr\}\mathbf{V}_{n}
\Biggr\rrVert _{2}^{2q}.
\end{eqnarray*}
An argument similar to that used to prove
\eqref{T1} implies
$E|(B3)|^{2q}=\mathrm{O}(k^{-q}(n-k)^{q+2qd})$.
In addition, by (\ref{4.11}), (\ref{4.12}), (\ref{resu1}) and $\sum_{t=1}^{n}\mathsf{q}_{t,i}^{2}=1$ for
$i=1,2,\ldots,r$, one obtains
\begin{eqnarray*}
E\bigl|(B4)\bigr|^{2q} \leq r^{2q-1}\sum_{i=1}^{r}
\bigl[E\bigl(v_{i}^{4q}\bigr) \bigr]^{1/2} \Biggl[E
\Biggl(\sum_{t=k}^{n-1}\mathsf{q}_{t+1-j,i}
\tau_{k,t} \Biggr)^{4q} \Biggr]^{1/2} \leq
Cn^{4qd},
\end{eqnarray*}
$E|(B5)|^{2q} \leq Cn^{4qd}$,
$E|(B6)|^{2q} \leq Cn^{4qd}$,
and $E(B7) \leq Ck^{q}n^{2qd}/n^{2q}$.
With the help of these moment inequalities, (\ref{lemma3.1-2}) and (\ref
{lemma4.2-1})--(\ref{lemma4.2-2}),
the proof of (i) can be completed in the same way as the proof of \eqref{T3}.
Moreover, by modifying
the proofs of \eqref{T4} and \eqref{T6} accordingly,
we can establish (ii) and (iii).
The details, however, are not presented here.
\end{pf}

%
\begin{prop}\label{proposition4.1}
Assume the same assumptions as in Proposition~\ref{proposition3.2}.
Suppose (\ref{conditionP}).
Then for any $\theta>1/q$,
\begin{longlist}[(ii)]
\item[(i)]
$\|\tilde{\mathbf{T}}_{n}(K_{n})-\mathbf{T}_{n}(K_{n})\|
_{2}^{2}=\mathrm{O}_{p}((K^{1+\theta}_{n}/n^{1-2d})+(K_{n}^{2+\theta}/n))$.
\item[(ii)]
\[
\bigl\|\tilde{\mathbf{D}}{}^{-1}_{n}(K_{n})-\mathbf
{D}_{n}^{-1}(K_{n})\bigr\|_{2}= \cases{
\mathrm{O}_{p}\bigl(n^{-1/2}K_{n}^{\theta}\bigr),&\quad$
\mbox{if } 0<d<1/4$;\vspace*{2pt}
\cr
\mathrm{O}_{p}\bigl((\log
n/n)^{1/2}K_{n}^{\theta}\bigr),&\quad$\mbox{if } d=1/4$;
\vspace *{2pt}
\cr
\mathrm{O}_{p}\bigl(n^{-1+2d}K_{n}^{\theta}
\bigr),&\quad$\mbox{if } 1/4<d<1/2$. }
\]
\end{longlist}
\end{prop}
\begin{pf}
In view of the proof of Proposition~\ref{proposition3.1},
(i) follows directly from (i) and (ii) of Lemma~\ref{lemma4.1}.
To show (ii),
note first that (\ref{proposition3.2-1})
still holds with $\hat{\sigma}_{0}^{2}$
replaced by $\tilde{\sigma}_{0}^{2}$.
This, in conjunction with (i)--(iii) of Lemma~\ref{lemma4.1}
and the argument used in the proof of Proposition~\ref{proposition3.2},
yields (ii).
\end{pf}
%
We are now in a position to introduce Theorem~\ref{theorem4.1}.
%
\begin{thmm}\label{theorem4.1}
Consider the regression model \eqref{regression
model}. With the same assumptions as in Proposition~\ref{proposition3.2},
suppose that \eqref{Kn} holds for some $\theta>1/q$.
Then
%
\begin{eqnarray} \label{thmm4.1}
&&\bigl\|\tilde{\bolds{\Omega }}{}^{-1}_{n}(K_{n})-
\bolds{\Omega}_{n}^{-1}\bigr\|_{2}
\nonumber
\\[-8pt]
\\[-8pt]
\nonumber
&&\quad =\mathrm{O}_{p} \biggl( \biggl(\frac{\log n \log
K_{n}}{K^{d}_{n}} \biggr)^{1/2} +
\biggl(\frac{K_{n}^{1+\theta} \log K_{n}}{n^{1-2d}}+ \frac{K_{n}^{2+\theta} \log K_{n}}{n} \biggr)^{1/2} \biggr)
=\mathrm{o}_{p}(1)
\end{eqnarray}
and
%
\begin{equation}
\bigl\|\tilde{\bolds{\Omega}}{}^{-1}_{n}(K_{n})
\bigr\|_{2}=\mathrm{O}_{p}(1). \label{thmm4.2}
\end{equation}
\end{thmm}
\begin{pf}
In view of the proof of Theorem~\ref{theorem3.1},
(\ref{thmm4.1}) and (\ref{thmm4.2})
are immediate consequences of
Proposition~\ref{proposition4.1} and Theorem~\ref{theorem2.2}.
\end{pf}
%
%
\begin{rem}\label{remark4.1}
Since no assumptions are imposed on the design matrix $\mathbf{X}_{n}$,
one of the most intriguing implications of Theorem~\ref{theorem4.1}
is that $\bolds{\Omega}_{n}^{-1}$ can be consistently estimated by
$\tilde{\bolds{\Omega}}{}^{-1}_{n}(K_{n})$ even when
$\mathbf{X}_{n}$ is singular.\vspace*{1pt}
Moreover, according to \eqref{thmm4.1}
and \eqref{thmm3.1}, it is interesting to point out that
$\tilde{\bolds{\Omega}}{}^{-1}_{n}(K_{n})$
and $\hat{\bolds{\Omega}}{}^{-1}_{n}(K_{n})$
share the same rate of convergence.
\end{rem}
%
Next, we consider the problem of estimating
$\mathbf{a}(n)$ under model \eqref{regression model}.
Recall Yule--Walker equations $\mathbf{a}(n)=\bolds{\Omega
}{}^{-1}_{n}\bolds{\gamma}_{n}$,
where $\bolds{\gamma}_{n}=(\gamma_{1},\ldots,\gamma_{n})'$.
A truncated version of
$\bolds{\Omega}_{n}^{-1}\bolds{\gamma}_{n}$
is given by
$\check{\mathbf{a}}(n)=\bolds{\Omega}_{n}^{-1}(K_{n})\check
{\bolds{\gamma}}_{n}$,
where $\check{\bolds{\gamma}}_{n}=(\gamma_{1},\ldots,\gamma
_{K_{n}},0,\ldots,0)'$
is an $n$-dimensional vector.
A natural estimate of
$\check{\mathbf{a}}(n)$
is
$\mathbf{a}^{*}(n)=\tilde{\bolds{\Omega}}{}^{-1}_{n}(K_{n})\tilde
{\bolds{\gamma}}_{n}$,
where $\tilde{\bolds{\gamma}}_{n}=(\tilde{\gamma}_{1},\ldots,\tilde
{\gamma}_{K_{n}},0,\ldots,0)'$
is an $n$-dimensional vector with
$\tilde{\gamma}_{j}$ denoting the $(1,j+1)$th entry of $\tilde
{\bolds{\Gamma}}_{K_{n}+1,n}$.
We shall show that when $K_{n}$ is suitably chosen,
$\mathbf{a}^{*}(n)$ is a consistent estimate of
$\mathbf{a}(n)$.
%
\begin{cor}\label{corollary4.1}
Assume the same assumptions as in Theorem~\ref{theorem4.1}. Suppose
that \eqref{Kn} holds and
$\frac{K_{n}^{1+2d/3}}{n^{1-2d}}=\mathrm{o}(1)$.
Then for any $\theta>1/q$,
\begin{eqnarray*}
\bigl\|\mathbf{a}^{*}(n)-\mathbf{a}(n)\bigr\|_{2} &=& \cases{
\displaystyle \mathrm{O}_{p}\biggl(\biggl(\frac{1}{K_{n}^{1-2d}}+\frac{\log n}{K_{n}^{d}}+
\frac
{K_{n}^{1+\theta}}{n^{1-2d}}+\frac{K_{n}^{2+\theta}}{n}\biggr)^{1/2}\biggr),&\quad$\mbox {if
}0<d\leq1/4$,\vspace*{2pt}
\cr
\displaystyle \mathrm{O}_{p}\biggl(\biggl(\frac{1}{K_{n}^{1-2d}}+
\frac{\log n}{K_{n}^{d}}+\frac
{K_{n}^{1+\theta}}{n^{1-2d}}+\frac
{K_{n}^{3+2d}}{n^{3-6d}}\biggr)^{1/2}
\biggr),&\quad$\mbox{if }1/4<d<1/2$}
\\
&=&\mathrm{o}_{p}(1).
\end{eqnarray*}
\end{cor}
\begin{pf}
Note first that
%
\begin{eqnarray}
\bigl\|\mathbf{a}^{*}(n)-\mathbf{a}(n)\bigr\|_{2}&\leq& \bigl\|\check{
\mathbf{a}}(n)-\mathbf{a}(n)\bigr\|_{2}+\bigl\|\mathbf{a}^{*}(n)-
\check {\mathbf{a}}(n)\bigr\|_{2},\label{app1-2}
\\
\bigl\|\check{\mathbf{a}}(n)-\mathbf{a}(n)\bigr\|_{2}&\leq& \bigl\|\bolds{
\Omega}_{n}^{-1}(\check{\bolds{\gamma }}_{n}-
\bolds{\gamma}_{n})\bigr\|_{2}+\bigl\|\bigl(\bolds{\Omega
}_{n}^{-1}(K_{n})-\bolds{
\Omega}_{n}^{-1}\bigr)\check{\bolds{\gamma
}}_{n}\bigr\|_{2},\label{app1-3}
\\
\bigl\|\mathbf{a}^{*}(n)-\check{\mathbf{a}}(n)\bigr\|_{2}&\leq&\bigl \|
\tilde{\bolds{\Omega}}{}^{-1}_{n}(K_{n}) (
\tilde{\bolds{\gamma }}_{n}-\check{\bolds{
\gamma}}_{n})\bigr\|_{2}+ \bigl\|\bigl(\tilde{\bolds{
\Omega}}{}^{-1}_{n}(K_{n})-\bolds{\Omega
}_{n}^{-1}(K_{n})\bigr)\check{\bolds{
\gamma}}_{n}\bigr\|_{2}.\label{app1-4}
\end{eqnarray}
Moreover,
%
\begin{eqnarray}
\bigl\|\bolds{\Omega}_{n}^{-1}(\check{\bolds{\gamma
}}_{n}-\bolds{\gamma}_{n})\bigr\|_{2} \leq\bigl\|
\mathbf{T}_{n}'\mathbf{D}_{n}^{-1}
\bigr\|_{2}\bigl\|\mathbf{T}_{n}(\check {\bolds{
\gamma}}_{n}-\bolds{\gamma}_{n})\bigr\|_{2} \leq
C\bigl\|\mathbf{T}_{n}(\check{\bolds{\gamma}}_{n}-
\bolds {\gamma}_{n})\bigr\|_{2}\label{app1-5}
\end{eqnarray}
and
%
\begin{eqnarray}\label{app1-6}
\bolds{\Omega}_{n}^{-1}-\bolds{
\Omega}_{n}^{-1}(K_{n}) &=&\bigl(
\mathbf{T}_{n}-\mathbf{T}_{n}(K_{n})
\bigr)'\mathbf{D}_{n}^{-1}\bigl(\bigl(\mathbf
{T}_{n}-\mathbf{T}_{n}(K_{n})\bigr)+
\mathbf{T}_{n}(K_{n})\bigr)
\nonumber
\\
&&{}+\mathbf{T}_{n}'(K_{n}) \bigl(
\mathbf{D}_{n}^{-1}-\mathbf {D}_{n}^{-1}(K_{n})
\bigr) \bigl(\bigl(\mathbf{T}_{n}-\mathbf{T}_{n}(K_{n})
\bigr)+\mathbf {T}_{n}(K_{n})\bigr)
\\
&&{}+\mathbf{T}_{n}'(K_{n})
\mathbf{D}_{n}^{-1}(K_{n}) \bigl(\mathbf
{T}_{n}-\mathbf{T}_{n}(K_{n})
\bigr).\nonumber
\end{eqnarray}
By \eqref{model-2}, \eqref{acf}, \eqref{yuleCon1}, \eqref{yuleCon2}, it
follows that
$\|\mathbf{T}_{n}(\check{\bolds{\gamma}}_{n}-\bolds{\gamma
}_{n})\|_{2}=\mathrm{O}(K_{n}^{-1/2+d})$,
$\|\mathbf{T}_{n}(K_{n})\check{\bolds{\gamma}}_{n}\|_{2}=\mathrm{O}(1)$, and
$\|(\mathbf{T}_{n}-\mathbf{T}_{n}(K_{n}))\check{\bolds{\gamma
}}_{n}\|_{2}=\mathrm{O}(K_{n}^{-1/2+d})$.
These bounds, together with \eqref{app1-3}, \eqref{app1-5} and \eqref
{app1-6}, yield
%
\begin{equation}
\bigl\|\check{\mathbf{a}}(n)-\mathbf{a}(n)\bigr\| _{2}=\mathrm{O}\bigl(K_{n}^{-1/2+d}+
\bigl(K_{n}^{-d}\log n\bigr)^{1/2}
\bigr).\label{app1-7}
\end{equation}
By the first moment bound theorem of Findley and Wei \cite{r5},
Cauchy--Schwarz inequality, Proposition~\ref{proposition4.1}, \eqref
{lemma4.2-2}, Lemma~\ref{lemma4.2}(ii),
\eqref{thmm4.2} and \eqref{app1-4},
it can be shown that
%
\begin{eqnarray}
\bigl\|\mathbf{a}^{*}(n)-\check{\mathbf{a}}(n)\bigr\|_{2}= \cases{
\mathrm{O}_{p}\bigl(\bigl(n^{-1+2d}K_{n}^{1+\theta}+n^{-1}K_{n}^{2+\theta}
\bigr)^{1/2}\bigr),&\quad$\mbox {if } 0<d\leq1/4$,\vspace*{2pt}
\cr
\mathrm{O}_{p}\bigl(\bigl(n^{-1+2d}K_{n}^{1+\theta}+n^{-3+6d}K_{n}^{3+2d}
\bigr)^{1/2}\bigr),&\quad$\mbox {if } 1/4<d<1/2$, } \label{app1-8}
\end{eqnarray}
for any $\theta>1/q$.
Now, the desired conclusion follows from \eqref{app1-2}, \eqref{app1-7}
and \eqref{app1-8}.
\end{pf}
%
%
\begin{rem}\label{remark4.2}
When $u_{1}, \ldots, u_{n}$
are observable,
Wu and Pourahmadi \cite{r17} constructed an estimate,
$\breve{\bolds{\Omega}}{}^{-1}_{n,K_{n}}\breve{\gamma}_{n}$,
of $\mathbf{a}(n)$,
where
$\breve{\bolds{\Omega}}_{n, K_{n}}=(\hat{\gamma}_{i-j}\mathbf
{1}_{|i-j|\leq K_{n}})_{1\leq i, j\leq n}$
and $\breve{\bolds{\gamma}}_{n}=(\breve{\gamma}_{1},\ldots,\breve
{\gamma}_{n})'$ with
$\breve{\gamma}_{i}=\hat{\gamma}_{i}\mathbf{1}_{\{i\leq K_{n}\}}$.
By assuming $\sum_{j=1}^{\infty}|\gamma_{j}|<\infty$, they obtained
a convergence rate of
the proposed estimate
in terms of $K_{n}$, the moment restriction on $w_{t}$,
and $\sum_{j=K_{n}}^{\infty}|\gamma_{k}|$;
see Corollary~2 of Wu and Pourahmadi \cite{r17}.
However, their proof,
relying heavily on $\sum_{j=1}^{\infty}|\gamma_{j}|<\infty$,
is no longer applicable here.
\end{rem}
%
\subsection{The rate of convergence of the FGLSE}\label{section4-2}
In this section, we assume that
$\mathbf{X}_{n}$ is nonsingular, and hence
$\bolds{\beta}$ is uniquely defined.
We estimate $\bolds{\beta}$
using the FGLSE,
\[
\hat{\bolds{\beta}}_{\mathrm{FGLS}}= \bigl(\mathbf{X}_{n}'
\tilde{\bolds{\Omega}}{}^{-1}_{n}(K_{n})
\mathbf{X}_{n}\bigr)^{-1} \mathbf{X}_{n}'
\tilde{\bolds{\Omega}}{}^{-1}_{n}(K_{n})
\mathbf{y}_{n}.
\]
The main objective of this section is to investigate
the convergence rate of $\hat{\bolds{\beta}}_{\mathrm{FGLS}}$.
To simplify the exposition, we shall focus
on polynomial regression models and impose
the following conditions on $a_{i}$:
%
\begin{eqnarray}
a_{j} \sim C_{0}j^{-1-d} \quad\mbox{and} \quad\sum
_{j=0}^{\infty}a_{j}\mathrm{e}^{ij\lambda}=0
\qquad\mbox{if and only if } \lambda=0\label{aiCond},
\end{eqnarray}
where $a_{0}=-1$ and $C_{0} \neq0$.
As mentioned in Section~\ref{section2},
\eqref{aiCond} is fulfilled by
the FARIMA model defined in \eqref{farima}.
When $K_{n}$ diverges to infinity at a suitable rate,
we derive the rate of convergence of $ \hat{\bolds{\beta}}_{\mathrm{FGLS}}$
in the next corollary.
It is important to be aware that
our proof is not a direct application of Theorem~\ref{theorem4.1}.
Instead, it relies on a very careful analysis
of the joint effects between
the Cholesky factors and the regressors.

%
\begin{cor}\label{corollary4.2}
Consider the regression model \eqref{regression model} with
$x_{ti}=t^{i-1}$ for $i=1,\ldots,p$. Assume the same
assumptions as in Theorem~\ref{theorem4.1} with \eqref{model-2}
replaced by \eqref{aiCond}. Suppose
that \eqref{conditionP} holds and
$n^{-1+2d}K_{n}^{1+2d}+n^{-1}K_{n}^{2+2d}=\mathrm{o}(1)$. Then
\begin{longlist}[(ii)]
\item[(i)]
$\|\mathbf{L}_{n}(\bolds{\beta}_{\mathrm{FGLS}}-\bolds{\beta})\|_{2}=\mathrm{O}_{p}(1)$,
\item[(ii)]
$\|\mathbf{L}_{n}(\hat{\bolds{\beta}}_{\mathrm{FGLS}}-\bolds{\beta})\|
_{2}=\mathrm{O}_{p}(1)$,
\end{longlist}
where $\mathbf{L}_{n}=n^{-d}\operatorname{diag}(n^{1/2},n^{3/2},\ldots,n^{p-1/2})$
and $\bolds{\beta}_{\mathrm{FGLS}}$ is $\hat{\bolds{\beta}}_{\mathrm{FGLS}}$
with $\tilde{\bolds{\Omega}}{}^{-1}_{n}(K_{n})$ replaced by
$\bolds{\Omega}_{n}^{-1}(K_{n})$.
\end{cor}
\begin{pf}
We only prove Corollary~\ref{corollary4.2} for $p=2$.
The proof for $p \neq2$ is analogous.
We begin
by showing (i). Let
$\tilde{\mathbf{L}}_{n}=K_{n}^{-d}\operatorname{diag}(n^{1/2},n^{3/2})$. Then
straightforward calculations yield
%
\begin{eqnarray}\label{app2-1}
&&\bigl\|\mathbf{L}_{n}(\bolds{\beta}_{\mathrm{FGLS}}-\bolds{
\beta})\bigr\| _{2}
\nonumber
\\[-8pt]
\\[-8pt]
\nonumber
&&\quad\leq n^{-d}K_{n}^{d} \bigl\|
\tilde{\mathbf{L}}_{n}\bigl(\mathbf{X}_{n}'
\bolds{\Omega }_{n}^{-1}(K_{n})
\mathbf{X}_{n}\bigr)^{-1}\tilde{\mathbf{L}}_{n}
\bigr\|_{2} \bigl\|\tilde{\mathbf{L}}{}^{-1}_{n}
\mathbf{X}_{n}'\bolds{\Omega }_{n}^{-1}(K_{n})
\mathbf{u}_{n}\bigr\|_{2}.
\end{eqnarray}
Moreover, by \eqref{aiCond},
%
\begin{eqnarray}
n^{-d}K_{n}^{d}\bigl\|\tilde{\mathbf{L}}{}^{-1}_{n}
\mathbf{X}_{n}'\bolds {\Omega}{}^{-1}_{n}(K_{n})
\mathbf{u}_{n}\bigr\|_{2}&=&\mathrm{O}_{p}(1),\label{app2-2}
\\
\bigl\|\tilde{\mathbf{L}}_{n}\bigl(\mathbf{X}_{n}'
\bolds{\Omega }{}^{-1}_{n}(K_{n})
\mathbf{X}_{n}\bigr){}^{-1}\tilde{\mathbf{L}}_{n}
\bigr\|_{2} &\leq& \Biggl(n{}^{-1}\sum_{t=0}^{\kappa n}
\lambda_{\min}(A_{\lfloor\kappa
n\rfloor+t} + A_{n-t})
\Biggr)^{-1} = \mathrm{O}(1) \label{app2-222},
\end{eqnarray}
where $\lambda_{\min}(A)$ denotes the minimum eigenvalue of matrix $A$,
$0<\kappa<1$ and $A_{t}=\eta_{t}'\eta_{t}$,
with $\eta_{t}$ denoting the $t$th row of $n^{1/2}\mathbf
{D}_{n}^{-1/2}(K_{n})\mathbf{T}_{n}(K_{n})\mathbf{X}_{n}\tilde{\mathbf
{L}}{}^{-1}_{n}$.
Combining \eqref{app2-1}--\eqref{app2-222} yields (i).
To show (ii), note first that
%
\begin{eqnarray}
\bigl\|\mathbf{L}_{n}(\hat{\bolds{\beta}}_{\mathrm{FGLS}}-
\bolds{\beta})\bigr\| _{2}&\leq& \bigl\|\mathbf{L}_{n}(\bolds{
\beta}_{\mathrm{FGLS}}-\bolds{\beta})\bigr\|_{2}+\bigl\|
\mathbf{L}_{n}(\hat{\bolds{\beta}}_{\mathrm{FGLS}}-\bolds{
\beta }_{\mathrm{FGLS}})\bigr\|_{2},\label{app2-3}
\\
\bigl\|\mathbf{L}_{n}(\hat{\bolds{\beta}}_{\mathrm{FGLS}}-
\bolds{\beta }_{\mathrm{FGLS}})\bigr\|_{2}&\leq&\bigl\|(D1)\bigr\|_{2}+
\bigl\|(D2)\bigr\|_{2},\label{app2-4}
\end{eqnarray}
where $(D1)=\mathbf{L}_{n}((\mathbf{X}_{n}'\tilde{\bolds{\Omega
}}{}^{-1}_{n}(K_{n})\mathbf{X}_{n}){}^{-1}-
(\mathbf{X}_{n}'\bolds{\Omega}{}^{-1}_{n}(K_{n})\mathbf
{X}_{n}){}^{-1})\mathbf{X}_{n}'\bolds{\Omega}{}^{-1}_{n}(K_{n})\mathbf{u}_{n}$,
and $(D2)=\mathbf{L}_{n}(\mathbf{X}_{n}'\tilde{\bolds{\Omega
}}{}^{-1}_{n}(K_{n})\mathbf{X}_{n}){}^{-1}
\mathbf{X}_{n}'(\tilde{\bolds{\Omega}}{}^{-1}_{n}(K_{n})-\bolds
{\Omega}{}^{-1}_{n}(K_{n}))\mathbf{u}_{n}$.
In addition,
%
\begin{eqnarray}
\bigl\|(D5)\bigr\|_{2}\leq\bigl(\bigl\|(D5)\bigr\|_{2}+\bigl\|(D3)\bigr\|_{2}
\bigr)\bigl\|(D4)\bigr\|_{2}\bigl\|(D3)\bigr\| _{2}\label{app2-5}
\end{eqnarray}
where $(D3)=\tilde{\mathbf{L}}_{n}(\mathbf{X}_{n}'\bolds{\Omega}{}^{-1}_{n}
(K_{n})\mathbf{X}{}^{-1}_{n})\tilde{\mathbf{L}}_{n}$,
$(D4)=\tilde{\mathbf{L}}{}^{-1}_{n}\mathbf{X}_{n}'(\tilde{\bolds
{\Omega}}{}^{-1}_{n}(K_{n})-
\bolds{\Omega}{}^{-1}_{n}(K_{n}))\mathbf{X}_{n}\tilde{\mathbf{L}}{}^{-1}_{n}$,
and
$(D5)=\tilde{\mathbf{L}}_{n}((\mathbf{X}_{n}'\tilde{\bolds{\Omega
}}{}^{-1}_{n}(K_{n})\mathbf{X}{}^{-1}_{n})-
(\mathbf{X}_{n}'\bolds{\Omega}{}^{-1}_{n}(K_{n})\mathbf
{X}{}^{-1}_{n}))\tilde{\mathbf{L}}_{n}$.
By \eqref{aiCond} and some algebraic manipulations, one obtains
$\|(D3)\|_{2}=\mathrm{O}(1)$ and $\|(D4)\|_{2}=\mathrm{o}_{p}(1)$.
Thus, by \eqref{app2-5}, $\|(D5)\|_{2}=\mathrm{o}_{p}(1)$.
The bounds for $\|(D3)\|_{2}$ and $\|(D5)\|_{2}$, together with \eqref
{aiCond} and \eqref{app2-2}, imply
%
\begin{eqnarray}
\bigl\|(D1)\bigr\|_{2}&\leq& n^{-d}K_{n}^{d}
\bigl\|(D5)\bigr\|_{2}\bigl\| \tilde{\mathbf{L}}{}^{-1}_{n}
\mathbf{X}_{n}'\bolds{\Omega }_{}^{-1}{n}(K_{n})
\mathbf{u}_{n}\bigr\|_{2}=\mathrm{o}_{p}(1), \label{hi}
\\
 \label{hi2}\bigl\|(D2)\bigr\|_{2}&\leq& n^{-d}K_{n}^{d}
\bigl(\bigl\|(D5)\bigr\|_{2}+\bigl\|(D3)\bigr\|_{2}\bigr) \bigl\|\tilde{
\mathbf{L}}{}^{-1}_{n}\mathbf{X}_{n}'
\bigl(\tilde{\bolds{\Omega }}{}^{-1}_{n}(K_{n})-
\bolds{\Omega}{}^{-1}_{n}(K_{n})\bigr)
\mathbf{u}_{n}\bigr\|_{2}\qquad
\nonumber
\\[-8pt]
\\[-8pt]
\nonumber
&=&\mathrm{o}_{p}(1).
\end{eqnarray}
Now, the desired conclusion (ii) follows from
\eqref{hi}, \eqref{hi2}, \eqref{app2-3} and \eqref{app2-4} and (i).
\end{pf}
%
%
\begin{rem}\label{remark4.3}
Under assumptions similar to those of Corollary~\ref{corollary4.2},
Theorems 2.2 and 2.3 of Yajima \cite{r18} show that
the best linear unbiased estimate (BLUE)
$\hat{\bolds{\beta}}_{\mathrm{BLUE}}=
(\mathbf{X}_{n}'\bolds{\Omega}{}^{-1}_{n}\mathbf{X}_{n}){}^{-1}
\mathbf{X}_{n}' \bolds{\Omega}{}^{-1}_{n}\mathbf{y}_{n}$,
and the
LSE,
$\hat{\bolds{\beta}}_{\mathrm{LS}}=
(\mathbf{X}_{n}'\mathbf{X}_{n}){}^{-1}
\mathbf{X}_{n}'\mathbf{y}_{n}$,
of $\bolds{\beta}$ have the same rate of convergence, and this
rate is, in turn, the same as that of $\hat{\bolds{\beta}}_{\mathrm{FGLS}}$.
\end{rem}

We close this section with a subtle example showing
that the convergence rate of $\hat{\bolds{\beta}}_{\mathrm{FGLS}}$
is faster than that of $\hat{\bolds{\beta}}_{\mathrm{LS}}$,
but slower than that of $\hat{\bolds{\beta}}_{\mathrm{BLUE}}$.
Consider model \eqref{regression model}, with $p=1$,
$x_{t1}=1+\cos(\theta t)$, and $\theta\neq0$. Assume the same
assumptions as in Corollary~\ref{corollary4.2}.
Then, by an argument similar to that used in the proof of
Corollary~\ref{corollary4.2}, it can be shown that
the rate of convergence of $\hat{\bolds{\beta}}_{\mathrm{FGLS}}$
is $n^{-1/2+d}K_{n}^{-d}$.
On the other hand,
Theorems 2.1 and 2.2 and Example~2.1(ii) of Yajima \cite{r19} yield that
the convergence rates of $\hat{\bolds{\beta}}_{\mathrm{BLUE}}$
and $\hat{\bolds{\beta}}_{\mathrm{LS}}$ are $n^{-1/2}$ and
$n^{-1/2+d}$, respectively.
This example gives a warning that
the convergence rate of $\hat{\bolds{\beta}}_{\mathrm{BLUE}}$
is not necessarily maintained by its feasible counterpart,
even if the consistency of $\tilde{\bolds{\Omega}}{}^{-1}_{n}(K_{n})$
holds true.

\section{Simulation study}\label{section5}
In Section~\ref{section5-1}, we introduce a data-driven method for
choosing the banding parameter $K_{n}$.
With this $K_{n}$,
we demonstrate
the finite sample performance of the inverse autocovariance estimator
proposed in Section~\ref{section3} under FARIMA($p,d,q$) processes,
and that
proposed in Section~\ref{section4} under
polynomial regression models
with I($d$) errors.
The details are given in
Sections~\ref{section5-2} and \ref{section5-3}, respectively.

\subsection{Selection of $K_{n}$}\label{section5-1}
Our approach for choosing $K_{n}$ is based on the idea of subsampling
and risk-minimization (SAR)
introduced by Bickel and Levina \cite{r3} and Wu and Pourahmadi \cite{r17}.
We first split the time series data $\{u_{i}\}_{i=1}^{n}$ into $\lfloor
n/b \rfloor$ nonoverlapping
subseries $\{u_{j}\}_{j=(v-1)b+1}^{vb}$ of equal length $b$,
where $b$ is a prescribed integer and $v=1,2,\ldots,\lfloor n/b \rfloor
$ with
$\lfloor a \rfloor$ denoting the largest integer $\leq a$.
Let $1\leq L<H<b$ be another prescribed integers.
For a given banding parameter $L\leq k<H$, let $\hat{\bolds{\Omega
}}{}^{-1}_{H,k,v}$ represent our
inverse autocovariance matrix estimator of
$\bolds{\Omega}_{H}{}^{-1}$
based on the
$v$th subseries $\{u_{j}\}_{j=(v-1)b+1}^{vb}$.
Define the average risk
\[
\hat{R}^{(O)}(k)=\frac{1}{\lfloor n/b \rfloor}\sum_{v=1}^{\lfloor n/b
\rfloor}
\bigl\llVert \hat{\bolds{\Omega}}{}^{-1}_{H,k,v}-\bolds{
\Omega }{}^{-1}_{H}\bigr\rrVert _{2}.
\]
Our goal is to find a banding parameter such that $\hat{R}^{(O)}(k)$ is
minimized.
However, since $\bolds{\Omega}{}^{-1}_{H}$ is unknown,
we use $\hat{\bolds{\Gamma}}{}^{-1}_{H,n}$, the $H$-dimensional
inverse sample autocovariance matrix, as its surrogate,
and replace
$\hat{R}^{(O)}(k)$ by
\[
\hat{R}(k)=\frac{1}{\lfloor n/b \rfloor}\sum_{v=1}^{\lfloor n/b
\rfloor}
\bigl\llVert \hat{\bolds{\Omega}}{}^{-1}_{H,k,v}-\hat{
\bolds {\Gamma}}{}^{-1}_{H,n}\bigr\rrVert _{2},
\]
noting that when $H \ll n$,
$\hat{\bolds{\Gamma}}{}^{-1}_{H,n}$
is a consistent estimator of $\bolds{\Omega}{}^{-1}_{H}$.
Now the banding parameter $K_{n}$ is chosen
to minimize $\hat{R}(k)$ over the interval $[L
,H)$.
In our simulation study, $b$ is set to $\lfloor n/5\rfloor$.
In addition, inspired by Theorem~\ref{theorem3.1},
we choose $L=\lfloor\log n \rfloor$
and $H=\lceil n^{0.4}\rceil$,
where $\lceil a \rceil$ denotes
the smallest integer $\geq a$.
The banding parameter for the detrended time series is also chosen in
the same manner.

\subsection{Finite sample performance of \texorpdfstring{$\hat{\bolds{\Omega}}{}^{-1}_{n}(K_{n})$}{hat{Omega}{-1}{n}(K{n})}}\label{section5-2}

We explore the finite sample performance of
$\hat{\bolds{\Omega}}{}^{-1}_{n}(K_{n})$, with $K_{n}$ determined by
the SAR method, under the following four data generating processes
(DGPs):
\begin{eqnarray*}
&&\mbox{DGP 1: }(1-B)^{d}u_{t}=w_{t};\qquad \mbox{DGP
2: }(1-0.7B) (1-B)^{d}u_{t}=w_{t};
\\
&&\mbox{DGP 3: }(1-B)^{d}u_{t}=(1-0.4B)w_{t}; \qquad
\mbox{DGP 4: }(1+0.4B) (1-B)^{d}u_{t}=(1-0.3B)w_{t},
\end{eqnarray*}
where
the $w_{t}$'s are i.i.d.
$N(0, 1)$ innovations.
To improve the speed and accuracy, we adopt
the method of Wu, Michailidis and Zhang \cite{r-1}
to generate the long memory data $\{u_{1}, \ldots, u_{n}\}$.
The performance of $\hat{\bolds{\Omega}}{}^{-1}_{n}(K_{n})$
is evaluated by $\hat{l}_{2}(d)$, the average value of
$\|\hat{\bolds{\Omega}}{}^{-1}_{n}(K_{n})-\bolds{\Omega
}{}^{-1}_{n}\|_{2}$ over 1000 replications,
with $n = 250,500,1000,2000,4000$.
The results are summarized in Table~\ref{App:Table1}.
Note first that for each combination of $d$ and DGP,
$\hat{l}_{2}(d)$ shows an obvious downward trend as $n$ increases.
Moreover, when $n=4000$, all $\hat{l}_{2}(d)$ are less than 0.65
except for $d=0.1$ and DGP${}={}$DGP 3 or DGP 4.
In the latter two cases,
$\hat{l}_{2}(d)$, lying between 0.93 and 1.34, are still reasonably small.
These findings suggest that
$\hat{\bolds{\Omega}}{}^{-1}_{n}(K_{n})$
is a reliable estimate of $\bolds{\Omega}{}^{-1}_{n}$,
particularly when $n$ is large enough.

%
\begin{table}
\caption{The values of $\hat{l}_{2}(d)$ under DGPs 1--4}\label{App:Table1}
\begin{tabular*}{\textwidth}{@{\extracolsep{\fill}}lccccccccccc@{}}
\hline
\multicolumn{1}{@{}l}{$n \setminus d$}
& \multicolumn{1}{l}{$0.01$} & \multicolumn{1}{l}{$0.1$} & \multicolumn
{1}{l}{$0.25$} & \multicolumn{1}{l}{$0.4$} & \multicolumn{1}{l}{$0.49$}
&& \multicolumn{1}{l}{$0.01$} & \multicolumn{1}{l}{$0.1$} &
\multicolumn{1}{l}{$0.25$} & \multicolumn{1}{l}{$0.4$} & \multicolumn
{1}{c}{$0.49$}
\\
\hline
&\multicolumn{5}{c}{DGP 1}&&\multicolumn{5}{c}{DGP 2}\\
\phantom{0}$250$ & 0.501 & 0.546 & 0.603 & 0.699 & 0.758 && 0.936 & 1.040 & 1.250
& 1.512 & 1.676 \\
\phantom{0}$500$ & 0.389 & 0.443 & 0.455 & 0.527 & 0.595 && 0.759 & 0.837 & 0.981
& 1.192 & 1.309 \\
$1000$ & 0.276 & 0.366 & 0.335 & 0.396 & 0.444 && 0.537 & 0.595 &
0.734 & 0.867 & 0.977 \\
$2000$ & 0.217 & 0.344 & 0.274 & 0.334 & 0.367 && 0.441 & 0.498 &
0.597 & 0.732 & 0.814 \\
$4000$ & 0.173 & 0.344 & 0.216 & 0.257 & 0.298 && 0.345 & 0.389 &
0.481 & 0.573 & 0.647 \\[3pt]
&\multicolumn{5}{c}{DGP 3}
&&\multicolumn{5}{c}{DGP 4}\\
\phantom{0}$250$ & 0.767 & 1.007 & 0.775 & 0.642 & 0.660 && 1.141 & 1.495 & 1.129
& 0.836 & 0.839 \\
\phantom{0}$500$ & 0.642 & 0.952 & 0.652 & 0.514 & 0.529 && 0.923 & 1.373 & 0.942
& 0.725 & 0.688 \\
$1000$ & 0.512 & 0.953 & 0.579 & 0.420 & 0.443 && 0.724 & 1.366 &
0.839 & 0.604 & 0.594 \\
$2000$ & 0.435 & 0.928 & 0.495 & 0.358 & 0.376 && 0.625 & 1.339 &
0.714 & 0.518 & 0.497 \\
$4000$ & 0.373 & 0.931 & 0.430 & 0.299 & 0.320 && 0.550 & 1.337 &
0.614 & 0.434 & 0.416 \\
\hline
\end{tabular*}
\end{table}

On the other hand,
the decreasing rate of $\hat{l}_{2}(d)$ apparently changes over $d$ and DGP.
To provide a better understanding of this phenomenon,
we first consider the fastest possible convergence rate
that can be derived from Theorem~\ref{theorem3.1}:
%
\begin{eqnarray}\label{T888}
&&\bigl\|\hat{\bolds{\Omega }}{}^{-1}_{n}\bigl(K^{*}_{n}
\bigr)-\bolds{\Omega}{}^{-1}_{n}\bigr\|_{2}
\nonumber
\\[-8pt]
\\[-8pt]
\nonumber
&&\quad= \cases{
 \mathrm{O}_{p}\bigl(n^{-{d}/{(4+2d+2\theta)}}
(\log n)^{{(4+d+2\theta)
}/{(4+2d+2\theta)}}\bigr),&\quad$\mbox{if }
0<d\leq\tilde{d}$,\vspace*{2pt}
\cr
 \mathrm{O}_{p}\bigl(n^{-{d(1-2d)}/{(2+2d+2\theta)}}(\log
n)^{{(2+d+2\theta)
}/{(2+2d+2\theta)}}\bigr),&\quad $\mbox{if } \tilde{d}<d<1/2$, }
\end{eqnarray}
where
$\tilde{d}=\{(3+2\theta)/2+\theta^{2}/4\}^{1/2}-(1+\theta/2)$, and
\[
K^{*}_{n}= \cases{ (n\log n)^{1/(2+d+\theta)},&\quad$\mbox{if }
0<d\leq\tilde{d}$,\vspace *{2pt}
\cr
(\log n)^{1/(1+d+\theta)}n^{(1-2d)/(1+\theta+d)},&\quad $
\mbox{if } \tilde{d}<d<1/2$. } 
\]
Because $w_{t}$'s are normally distributed, in view of Theorem~\ref{theorem3.1}, $\theta$ can be any positive number, and hence
$\tilde{d}$ is arbitrarily close to $\sqrt{1.5}-1$ (which, rounded
to the nearest thousandth, is 0.225).
We then measure the relative performance of
$\|\hat{\bolds{\Omega}}{}^{-1}_{n}(K_{n})-\bolds{\Omega
}{}^{-1}_{n}\|_{2}$
and
$\|\hat{\bolds{\Omega}}{}^{-1}_{n}(K^{*}_{n})-\bolds{\Omega
}{}^{-1}_{n}\|_{2}$
using the ratio $\hat{l}_{2}(d)/\operatorname{ OP}(d)$, where
%
\begin{eqnarray}
\operatorname{ OP}(d)= \cases{ 0.05n^{-{d}/{(4+2d)}}
(\log n)^{{(4+d)}/{(4+2d)}},&\quad$\mbox{if }
0<d \leq0.225$,\vspace*{2pt}
\cr
0.05n^{-{d(1-2d)}/{(2+2d)}}(\log n)^{{(2+d)}/{(2+2d)}},&\quad$
\mbox{if } 0.225<d<1/2$, } \label{T88888}
\end{eqnarray}
which is obtained from the bound in (\ref{T888}) with $\theta$ set to 0
and constants set to
0.05.
The values of
$\hat{l}_{2}(d)/\operatorname{ OP}(d)$ under DGPs 1--4 are summarized in Table~\ref{Table2}.
Note that while the exact constants are not reported in \eqref{T888},
setting them to 0.05 helps us
to better interpret some numerical results in Table~\ref{App:Table1}
through Table~\ref{Table2}.

\begin{table}
\caption{The values of $\hat{l}_{2}(d)/\operatorname{ OP}(d)$ under the DGPs
1--4}\label{Table2}
\begin{tabular*}{\textwidth}{@{\extracolsep{\fill}}lccccccccccc@{}}
\hline
\multicolumn{1}{@{}l}{$n \setminus d$}
& \multicolumn{1}{l}{$0.01$} & \multicolumn{1}{l}{$0.1$} & \multicolumn
{1}{l}{$0.25$} & \multicolumn{1}{l}{$0.4$} & \multicolumn{1}{l}{$0.49$}
&& \multicolumn{1}{l}{$0.01$} & \multicolumn{1}{l}{$0.1$} &
\multicolumn{1}{l}{$0.25$} & \multicolumn{1}{l}{$0.4$} & \multicolumn
{1}{c}{$0.49$}
\\
\hline
&\multicolumn{5}{c}{DGP 1}
&&\multicolumn{5}{c}{DGP 2}\\
\phantom{0}$250$ & 1.849 & 2.349 & 3.414 & 3.782 & 3.704 && 3.453 & 4.476 & 7.081
& 8.187 & 8.187 \\
\phantom{0}$500$ & 1.278 & 1.728 & 2.399 & 2.629 & 2.640 && 2.493 & 3.264 & 5.172
& 5.950 & 5.804 \\
$1000$ & 0.817 & 1.309 & 1.661 & 1.842 & 1.805 && 1.589 & 2.127 &
3.644 & 4.030 & 3.976 \\
$2000$ & 0.585 & 1.138 & 1.291 & 1.460 & 1.380 && 1.188 & 1.647 &
2.814 & 3.198 & 3.064 \\
$4000$ & 0.428 & 1.062 & 0.973 & 1.062 & 1.045 && 0.854 & 1.200 &
2.171 & 2.370 & 2.272 \\[3pt]
&\multicolumn{5}{c}{DGP 3}
&&\multicolumn{5}{c}{DGP 4}\\
\phantom{0}$250$ & 2.829 & 4.333 & 4.388 & 3.476 & 3.226 && 4.207 & 6.433 & 6.392
& 4.528 & 4.096 \\
\phantom{0}$500$ & 2.108 & 3.710 & 3.434 & 2.565 & 2.345 && 3.032 & 5.352 & 4.965
& 3.618 & 3.052 \\
$1000$ & 1.514 & 3.405 & 2.873 & 1.951 & 1.802 && 2.143 & 4.880 &
4.162 & 2.807 & 2.418 \\
$2000$ & 1.172 & 3.072 & 2.333 & 1.564 & 1.415 && 1.685 & 4.431 &
3.366 & 2.263 & 1.873 \\
$4000$ & 0.922 & 2.875 & 1.939 & 1.237 & 1.122 && 1.361 & 4.130 &
2.772 & 1.795 & 1.461 \\
\hline
\end{tabular*}
\end{table}

For $n \geq1000$,
all values of $\hat{l}_{2}(d)/\operatorname{ OP}(d)$
fall in a reasonable range of $(0.4, 5.0)$,
suggesting that the rate of convergence
of $\|\hat{\bolds{\Omega}}{}^{-1}_{n}(K_{n})-\bolds{\Omega
}{}^{-1}_{n}\|_{2}$
is comparable to the optimal rate
obtained from Theorem~\ref{theorem3.1}.
Moreover, the asymptotic behaviors of $\hat{l}_{2}(d)$
can be well explained by $\operatorname{ OP}(d)$
when DGP${}={}$DGP 1 and $d\geq0.1$.
In particular, when $n=4000$, the rankings of $\{\hat{l}_{2}(0.1),
\hat{l}_{2}(0.25), \hat{l}_{2}(0.4), \hat{l}_{2}(0.49)\}$ coincide
exactly with those of $\{\operatorname{ OP}(0.1), \operatorname{ OP}(0.25),
\operatorname{
OP}(0.4),  \operatorname{ OP}(0.49)\}$, and
$\operatorname{OP}(0.25)=\min_{d \in\{0.1, 0.25, 0.4, 0.49\}}\operatorname{OP}(d)$. This
gives reasons for explaining why $d=0.25$ often provides better
results than $d=0.1, 0.4$ or 0.49.
The behavior of $\hat{l}_{2}(0.01)$, however,
is apparently inconsistent with that of $\operatorname{ OP}(0.01)$. Specifically,
for $n\geq250$, $\hat{l}_{2}(0.01)< \min_{d\in\{0.1, 0.25, 0.4, 0.49\}
}\hat{l}_{2}(d)$,
whereas $\operatorname{OP}(0.01)> \max_{d \in\{0.1, 0.25, 0.4, 0.49\}}\operatorname{OP}(d)$.
One possible explanation of this discrepancy is that
when $d$ is extremely small,
the constant associated with the convergence rate
of $\|\hat{\bolds{\Omega}}{}^{-1}_{n}(K_{n})-\bolds{\Omega
}{}^{-1}_{n}\|_{2}$
can also be very small,
and the constant, 0.05, assigned to $\operatorname{ OP}(d)$
fails to do a good job in this extremal case.

It is relatively difficult to understand
the behaviors of $\hat{l}_{2}(d)$ through $\operatorname{ OP}(d)$
when short-memory AR or MA components are added into the I($d$) model.
However, using the $\hat{l}_{2}(d)$ in DGP 1 as the basis for comparison,
it seems fair to comment that
the AR component tends to increase $\hat{l}_{2}(d)$ with $d\geq0.25$
and $d=0.01$,
whereas the MA component tends to
increase $\hat{l}_{2}(d)$ with $d \leq0.25$.
When both components are included,
the values of $\hat{l}_{2}(d)$
are uniformly larger than those in the I($d$) case.
We leave a further investigation of the impact of the AR and MA components
on the finite sample performance of
$\hat{\bolds{\Omega}}{}^{-1}_{n}(K_{n})$
as a future work.

In the following, we shall perform a sensitivity analysis of the SAR method
by perturbing the parameter $c$ in $cK_{n}$.
We define the sensitivity function
\[
\operatorname{ SF}(c)=\frac{\|\hat{\bolds{\Omega
}}{}^{-1}_{n}(cK_{n})-\bolds{\Omega}{}^{-1}_{n}\|_{2}-\|\hat
{\bolds{\Omega}}{}^{-1}_{n}(K_{n})
-\bolds{\Omega}{}^{-1}_{n}\|_{2}}{\|\hat{\bolds{\Omega}}{}^{-1}_{n}(K_{n})
-\bolds{\Omega}{}^{-1}_{n}\|_{2}}.
\]
For each $c=0.8, 1.2$, $d=0.1, 0.25, 0.45$, DGP${}={}$DGP 1--4, and
$n=250, 500, 1000, 2000, \break  4000$, we compute the average of $\operatorname{
SF}(c)$, denoted by $\overline{\operatorname{ SF}}(c)$, based on 1000 replications,
and the five-number summaries
of $\overline{\operatorname{ SF}}(c)$
for each $n$ are presented in Table~\ref{App:Table4}.
Table~\ref{App:Table4} shows that
the maximum values of $\overline{\operatorname{ SF}}(c)$
are all positive and decrease as $n$ increases.
In contrast, the minimum values of $\overline{\operatorname{ SF}}(c)$
are all negative and start to increase when $n \geq1000$.
When $n=4000$,
the maximum $\overline{\operatorname{ SF}}(c)$ and minimum $\overline{\operatorname{ SF}}(c)$
are 0.152 and $-0.194$, respectively,
yielding that the average of
$\|\hat{\bolds{\Omega}}{}^{-1}_{n}(cK_{n})-\bolds{\Omega
}{}^{-1}_{n}\|_{2}$
falls between
0.806--1.152 times the average of $\|\hat{\bolds{\Omega
}}{}^{-1}_{n}(K_{n})-\bolds{\Omega}{}^{-1}_{n}\|_{2}$,
for all $c$'s, $d$'s and DGPs under consideration.
Our analysis reveals
that a small perturbation of $K_{n}$
will not lead to a drastic change on estimation errors.
%
\begin{table}
\caption{5-number summaries of $\overline{\operatorname{ SF}}(c)$}\label{App:Table4}
\begin{tabular*}{\textwidth}{@{\extracolsep{\fill}}ld{2.3}d{2.3}d{2.3}cc@{}}
\hline
\multicolumn{1}{@{}l}{$n$}
& \multicolumn{1}{l}{Minimum} & \multicolumn{1}{l}{1st quartile} & \multicolumn{1}{l}{Median} & \multicolumn{1}{l}{3rd quartile} & \multicolumn{1}{@{}l}{Maximum}\\
\hline
\phantom{0}$250$ & -0.173 & -0.109 & 0.050 & 0.087 & 0.289 \\
\phantom{0}$500$ & -0.229 & -0.074 & 0.044 & 0.150 & 0.228 \\
$1000$ & -0.250 & -0.022 & 0.029 & 0.160 & 0.210 \\
$2000$ & -0.210 & -0.025 & -0.005 & 0.134 & 0.156 \\
$4000$ & -0.194 & -0.029 & -0.014 & 0.112 & 0.152 \\
\hline
\end{tabular*}
\end{table}

\subsection{Finite sample performance of \texorpdfstring{$\tilde{\bolds{\Omega}}{}^{-1}_{n}(K_{n})$}{tilde{Omega}{-1}{n}(K{n})}}\label{section5-3}

We consider three polynomial regression models:
\begin{longlist}[Model 3:]
\item[Model 1:] $y_{t}=1+u_{t},  t=1,2,\ldots,n$,

\item[Model 2:] $y_{t}=1+2t+u_{t},  t=1,2,\ldots,n$,

\item[Model 3:] $y_{t}=5+t+2t^{4}+u_{t},  t=1,2,\ldots,n$,
\end{longlist}
where
$u_{t}$'s are generated by DGP 1.
The performance of
$\tilde{\bolds{\Omega}}{}^{-1}_{n}(K_{n})$ (with $K_{n}$ determined
by the SAR method)
is investigated with polynomial degree known or unknown.
In the latter situation,
we perform best subset selection in the following fifth-order model,
\[
y_{t}=\beta_{0}+\beta_{1}t+
\beta_{2}t^{2}+\beta_{3}t^{3}+\beta
_{4}t^{4}+\beta_{5}t^{5}+u_{t},\qquad
t=1,2,\ldots,n,
\]
using the selection criterion,
%
\begin{equation}
L_{n}(M)=\log\hat{\sigma}_{n}^{2}(M)+\#M/
\log(n),\label{lossFun}
\end{equation}
suggested by Ing and Wei \cite{r11},
where
$\mathcal{M}=\{M\dvt M\subseteq\{1,t,t^{2},t^{3},t^{4},t^{5}\}\}$
and $\hat{\sigma}_{n}^{2}(M)$ is the
residual mean square error of model $M$.
Note that according to
Theorem~4.1 of Ing and Wei \cite{r11},
$L_{n}(M)$ is a consistent criterion
in regression models with long-memory errors.
The performance of
$\tilde{\bolds{\Omega}}{}^{-1}_{n}(K_{n})$
is evaluated by $\tilde{l}_{2}(d)$,
which is $\hat{l}_{2}(d)$ with
$u_{t}$'s replaced by the corresponding detrended series.
The values of $\tilde{l}_{2}(d)$ are documented
in Table~\ref{App:Table3}, in which
$d\in\{0.1, 0.25, 0.4\}$, $n \in\{250, 500, 1000, 2000, 4000\}$
and models are known or selected by $L_{n}(M)$.
Table~\ref{App:Table3} also reports
the correct selection frequencies (in 1000 simulations), which is
denoted by $\hat{q}_{i}(d)$
for model~$i$ and long-memory parameter $d$.

%
\begin{table}
\caption{$\tilde{l}_{2}(d)$ with (in parentheses) or without model
selection and $\hat{q}_{i}(d)$}%
\label{App:Table3}
\begin{tabular*}{\textwidth}{@{\extracolsep{\fill}}lccccccc@{}}
\hline
&\multicolumn{7}{c}{Model 1}\\
\hline
&\multicolumn{3}{l}{$\tilde{l}_{2}(d)$}&&\multicolumn{3}{l}{$\hat
{q}_{1}(d)$}\\[-6pt]
&\multicolumn{3}{l}{\hrulefill}&&\multicolumn{3}{l@{}}{\hrulefill}\\
\multicolumn{1}{@{}l}{$n \setminus d$}
& \multicolumn{1}{l}{$0.1$} & \multicolumn{1}{l}{$0.25$} & \multicolumn
{1}{l}{$0.45$}
&& \multicolumn{1}{l}{$0.1$} & \multicolumn{1}{l}{$0.25$} &
\multicolumn{1}{l}{$0.45$}\\
\hline
\phantom{0}$250$ & 0.566 (0.566) & 0.592 (0.593) & 0.702 (0.693) && 0.992 & 0.867 &
0.558\\
\phantom{0}$500$ & 0.461 (0.461) & 0.456 (0.454) & 0.535 (0.533) && 0.999 & 0.918 &
0.593\\
$1000$ & 0.385 (0.385) & 0.333 (0.332) & 0.404 (0.400) && 1.000 & 0.962 &
0.622\\
$2000$ & 0.358 (0.358) & 0.271 (0.271) & 0.334 (0.331) && 1.000 & 0.984 &
0.620\\
$4000$ & 0.353 (0.353) & 0.216 (0.216) & 0.267 (0.264) && 1.000 & 0.996 &
0.609\\[6pt]
&\multicolumn{7}{c}{Model 2}\\
\hline
&\multicolumn{3}{l}{$\tilde{l}_{2}(d)$}&&\multicolumn{3}{l}{$\hat
{q}_{2}(d)$}\\[-6pt]
&\multicolumn{3}{l}{\hrulefill}&&\multicolumn{3}{l@{}}{\hrulefill}\\
\multicolumn{1}{@{}l}{$n \setminus d$}
& \multicolumn{1}{l}{$0.1$} & \multicolumn{1}{l}{$0.25$} & \multicolumn
{1}{l}{$0.45$}
&& \multicolumn{1}{l}{$0.1$} & \multicolumn{1}{l}{$0.25$} &
\multicolumn{1}{l}{$0.45$}\\
\hline
\phantom{0}$250$ & 0.585 (0.574) & 0.602 (0.604) & 0.690 (0.691) && 0.688 & 0.487 &
0.455\\
\phantom{0}$500$ & 0.477 (0.470) & 0.457 (0.455) & 0.533 (0.532) && 0.839 & 0.552 &
0.454\\
$1000$ & 0.399 (0.397) & 0.336 (0.336) & 0.400 (0.398) && 0.947 & 0.657 &
0.475\\
$2000$ & 0.368 (0.368) & 0.273 (0.275) & 0.331 (0.331) && 0.995 & 0.735 &
0.444\\
$4000$ & 0.359 (0.359) & 0.218 (0.218) & 0.263 (0.264) && 1.000 & 0.809 &
0.458\\[6pt]
&\multicolumn{7}{c}{Model 3}\\
\hline
&\multicolumn{3}{l}{$\tilde{l}_{2}(d)$}&&\multicolumn{3}{l}{$\hat
{q}_{3}(d)$}\\[-6pt]
&\multicolumn{3}{l}{\hrulefill}&&\multicolumn{3}{l@{}}{\hrulefill}\\
\multicolumn{1}{@{}l}{$n \setminus d$}
& \multicolumn{1}{l}{$0.1$} & \multicolumn{1}{l}{$0.25$} & \multicolumn
{1}{l}{$0.45$}
&& \multicolumn{1}{l}{$0.1$} & \multicolumn{1}{l}{$0.25$} &
\multicolumn{1}{l}{$0.45$}\\
\hline
\phantom{0}$250$ & 0.611 (0.611) & 0.617 (0.618) & 0.687 (0.686) && 1.000 & 0.995 &
0.904\\
\phantom{0}$500$ & 0.493 (0.493) & 0.463 (0.463) & 0.532 (0.530) && 1.000 & 0.999 &
0.912\\
$1000$ & 0.411 (0.411) & 0.339 (0.339) & 0.396 (0.395) && 1.000 & 1.000 &
0.916\\
$2000$ & 0.376 (0.376) & 0.275 (0.275) & 0.329 (0.328) && 1.000 & 1.000 &
0.942\\
$4000$ & 0.318 (0.318) & 0.218 (0.218) & 0.262 (0.261) && 1.000 & 1.000 &
0.967\\
\hline
\end{tabular*}
\end{table}

All $\hat{q}_{3}(d)$'s are larger than 0.9.
However,
$\hat{q}_{1}(0.45)$ and $\hat{q}_{2}(0.45)$
only fall in the interval (0.44, 0.63)
and the intercept (constant time trend)
is often excluded by $L_{n}(M)$ in these cases.
In fact,
identifying the intercept is a notoriously challenging problem when
$d$ is large and the intercept parameter is not far enough away from 0.
Fortunately,
Table~\ref{App:Table3} shows that
the $\tilde{l}_{2}(d)$ values obtained with or without model selection procedure
are similar, even when $\hat{q}_{i}(d)$
is much smaller than 1.
This result may be due to the fact that
under models 1 and 2,
the performance of
$\tilde{\bolds{\Omega}}{}^{-1}_{n}(K_{n})$
is insensitive to misspecification of the intercept, provided
$d$ is large enough.
Another interesting finding is that
for each regression model considered in this section
and each $(n, d)$ combination,
the behavior of $\tilde{l}_{2}(d)$
coincides with that of $\hat{l}_{2}(d)$ with DGP${}={}$DGP 1.
Putting these characteristics together suggests that
$\tilde{\bolds{\Omega}}{}^{-1}_{n}(K_{n})$
is a reliable surrogate for
$\hat{\bolds{\Omega}}{}^{-1}_{n}(K_{n})$. This conclusion
is particularly relevant in situations where
the latter matrix becomes infeasible.

\section*{Acknowledgements}
The work of Ching-Kang Ing was supported in part by
the National Science Council of Taiwan under Grant NSC
97-2628-M-001-022-MY2 and Academia Sinica Investigator Award. The
research of Hai-Tang Chiou was supported by
Academia Sinica Investigator Award and NSC 102-2118-M-110-002-MY2 from
Taiwan's National Science Council.
The research of Meihui Guo was partly supported by Grant number
NSC 102-2118-M-110-002-MY2, from Taiwan's National Science Council.

%

%




\printhistory
\end{document}